\documentclass[oneside,onecolumn,12pt]{article}

\usepackage{geometry}
\usepackage[utf8]{inputenc}
\usepackage[T1]{fontenc}

\usepackage{amsmath,amssymb,amsfonts,amsthm}
\usepackage{mathrsfs}   
\usepackage{romannum}

\usepackage{graphicx}   
\usepackage{subcaption} 
\usepackage{float}
\usepackage{array}      

\theoremstyle{plain}
\newtheorem{theorem}{Theorem}[section]
\newtheorem{lemma}{Lemma}[section]

\newcommand{\norm}[1]{\left\| #1 \right\|}     

\usepackage{marvosym} 
\usepackage{fancyhdr}
\usepackage{hyperref}
\hypersetup{
    colorlinks=true,    
    linkcolor=blue,     
    citecolor=blue,     
    urlcolor=blue,      
    breaklinks=true     
}

\begin{document}

\pagenumbering{arabic}
\setcounter{page}{1}

\title{Strong convergence of a slowly damped primal-dual
dynamical system with general time-scaling, Tikhonov
regularization and implicit Hessian damping}

\author{
Yi-xuan Ma\textsuperscript{a},
Hong-lu Li\textsuperscript{a},
and Yi-bin Xiao\textsuperscript{a,*}
}

\date{November 1, 2026}

\maketitle

\begingroup
\renewcommand{\thefootnote}{}
\footnotetext{%
\textsuperscript{a}Department of Mathematics,
University of Electronic Science and Technology of China,
Chengdu, Sichuan, P.R. China.\\
\textsuperscript{*}Corresponding author.\\
E-mail addresses:
\texttt{mayixuan929@outlook.com} (Y. Ma);
\texttt{lihonglu\_7988@163.com} (H. Li);
\texttt{xiaoyb9999@hotmail.com} (Y. Xiao).
}
\endgroup

\begin{abstract}
This paper proposes a novel mixed-order inertial primal-dual dynamical system with general time-scaling and Tikhonov regularization coefficients, together with implicit Hessian damping, for solving convex optimization problems with linear equality constraints in Hilbert spaces. The proposed dynamical system incorporates slowly vanishing viscous damping and is developed based on a time-dependent Tikhonov-regularized Lagrangian. Under suitable assumptions, we establish fast convergence rates for the objective residual, the feasibility violation, and the primal-dual gap using a Lyapunov analysis. Moreover, we prove that the generated primal-dual trajectory converges strongly to the minimum-norm primal-dual solution. Numerical experiments illustrate the theoretical results and indicate improved convergence performance relative to the comparison dynamical systems.
\end{abstract}

\begin{center}
\end{center}

\section{Introduction}
Let $\mathcal{X}$ and $\mathcal{Y}$ be two real Hilbert spaces endowed with the inner product $\langle \cdot, \cdot \rangle$ and the induced norm $\| \cdot \|$. This paper studies the following convex optimization problem with linear equality constraints:
\begin{equation} \label{uncon}
\min_{x \in \mathcal{X}} f(x), \quad \text{s.t.} \ Ax = b, 
\end{equation}
where $f : \mathcal{X} \to \mathbb{R}$ is a continuously differentiable convex function, $A : \mathcal{X} \to \mathcal{Y}$ is a continuous linear operator, and $b \in \mathcal{Y}$.
The feasible set of problem \eqref{uncon} is denoted by \(
\mathbb{F} := \{ x \in \mathcal{X} : Ax = b \}.
\)
When $A = 0$ and $b = 0$, problem \eqref{uncon} reduces to the unconstrained convex optimization problem
\begin{equation}\label{without}
\min_{x \in \mathcal{X}} f(x).
\end{equation}
Problem \eqref{uncon} serves as a fundamental model in numerous interdisciplinary areas such as communication systems, machine learning, and image restoration. For related applications, we refer the readers to \cite{wei_modeling_2026,przybylo_noninvasive_2026,vijayan_cahmo_2026,ma_integrated_2026,zhang_dynamic_2026}.

Associated with problem \eqref{uncon}, the Lagrangian function $\mathcal{L} : \mathcal{X} \times \mathcal{Y} \to \mathbb{R}$ associated with problem \eqref{uncon} is defined by
\[
\mathcal{L}(x, \lambda) = f(x) + \langle \lambda, Ax - b \rangle.
\]
Let $\Omega$ denote the set of saddle points of $\mathcal{L}$, 
\[
\Omega:=\left\{(\tilde{x},\tilde{\lambda})\in \mathcal{X}\times \mathcal{Y}:
\mathcal{L}(\tilde{x},\lambda)
\leq
\mathcal{L}(\tilde{x},\tilde{\lambda})
\leq
\mathcal{L}(x,\tilde{\lambda}),
\quad
\forall (x,\lambda)\in \mathcal{X}\times\mathcal{Y}
\right\}.
\]
Every $(x^*,\lambda^*)\in\Omega$ satisfies the primal-dual
optimality conditions. Throughout this paper, we assume that $\Omega\neq\varnothing$.
\subsection{Related works}

\quad  \ Continuous time dynamical systems provide a useful framework for analyzing optimization methods. They make it possible to investigate the respective roles of inertia, damping, time scaling, and regularization at the trajectory level, where Lyapunov methods can be used to reveal the mechanisms underlying convergence and acceleration. Moreover, suitably designed continuous time models often provide valuable guidance for the construction and analysis of discrete optimization algorithms. Early inertial continuous-time models for optimization can be traced back to Polyak's heavy-ball dynamical system (HBF) \cite{polyak_introduction_1987,polyak_methods_1964}. A landmark development in the modern continuous-time analysis of accelerated optimization was made by Su, Boyd, and Candès \cite{su_differential_2016}, who derived the second-order differential equation as the continuous-time limit of Nesterov's accelerated gradient method. Their work established a direct connection between accelerated first-order algorithms and inertial differential equations and provided a dynamical interpretation of Nesterov acceleration, thereby stimulating extensive subsequent research on accelerated optimization dynamics. Further convergence results for (HBF) in the cases of differentiable, quadratic, convex, and non-convex objective functions can be found in \cite{polyak_introduction_1987,felipe_alvarez_minimizing_2000,begout_damped_2015}. Compared with (HBF), Cabot et al. \cite{cabot_long_2009} proposed a dynamical system of a more general form :
\[
\ddot{x}(t) + a(t) \dot{x}(t) + \nabla f(x(t)) = 0,
\]
where $a(t)$ is a positive, asymptotically vanishing damping coefficient. Its rate of convergence and perturbed counterpart can be found in \cite{sebbouh_convergence_2020,attouch_fast_2018,attouch_rate_2019}. However, inertial dynamical systems governed solely by viscous damping may exhibit pronounced oscillations. To alleviate this phenomenon, Alvarez et al. \cite{alvarez_second_order_2002} introduced a second-order gradient-like dynamical system with Hessian-driven damping. Subsequently, Attouch et al. \cite{attouch_fast_2016} combined Hessian-driven damping with asymptotically vanishing viscous damping. Beyond constant or prescribed power-law coefficients, allowing the damping and time scaling parameters to be general functions provides a more flexible and unified framework. In particular, suitably chosen time-scale factors may further accelerate the dynamics, and the flexibility of the parameter functions can also improve the transient behavior of the trajectories.  Attouch et al. \cite{attouch_first-order_2022} showed that the introduction of time-scale factors can further accelerate the associated first-order optimization methods. In constrained mixed-order dynamics, appropriate choices of the time-scaling coefficient may even lead to exponential convergence rates \cite{he_second-order_2022}. Motivated by these advantages, Attouch et al. \cite{attouch_first-order_2022,attouch_convergence_2023} considered inertial systems with general time-dependent viscous damping, Hessian damping, and scaling coefficients:

\begin{equation*}
\ddot{x}(t)+\gamma(t)\dot{x}(t)
+\beta(t)\nabla^2 f(x(t))\dot{x}(t)
+b(t)\nabla f(x(t))=0,
\end{equation*}
where $\gamma(t)$, $\beta(t)$, and $b(t)$ are damping and scaling parameters. Under suitable conditions, these studies established fast convergence estimates for the objective values and weak convergence of the trajectories. Nevertheless, explicit Hessian damping requires the objective function to be twice continuously differentiable, which restricts its applicability, and such dynamical systems do not naturally discretize to Nesterov-type inertial algorithms. To overcome the drawback, László \cite{laszlo_solving_2024} proposed the following inertial dynamical system with implicit Hessian damping. Under suitable parameter conditions, László \cite{laszlo_solving_2024} obtained fast convergence rates for the objective values and established both weak and strong convergence results for the trajectories. 
 
Although many inertial dynamical systems achieve fast convergence of objective values, strong convergence of trajectories is generally difficult to obtain in Hilbert spaces. In fact, Baillon's \cite{baillon_exemple_1978} counter example shows that trajectories generated by the continuous steepest descent flow may converge weakly but not strongly. Attouch and Czarnecki \cite{attouch_asymptotic_2002} investigated a time-dependent Tikhonov regularization of (HBF). They proved that each trajectory strongly converges to the minimum-norm minimizer of $\min_{x \in \mathcal{X}} \Phi(x)$ when $\Phi$ is convex and the solution set $\arg\min\Phi$ is nonempty. Building on this Tikhonov selection principle, Attouch et al. \cite{attouch_damped_2022} investigated the following regularized inertial dynamical system with asymptotically vanishing damping:
\begin{equation*}
\ddot{x}(t)
+
\frac{\delta}{t^{p/2}}\dot{x}(t)
+
\nabla f(x(t))
+
\frac{1}{t^p}x(t)=0,
\end{equation*}
where $0<p<2$. Under suitable assumptions, Attouch et al. proved strong convergence of the trajectories together with convergence rate estimates.

The above developments mainly concern unconstrained convex optimization. Since many optimization models arising in applications are subject to linear equality constraints, it is natural to extend inertial dynamical approaches to the constrained setting. In this case, the dynamics must simultaneously account for objective minimization and feasibility. Introducing the Lagrange multiplier associated with the constraint naturally leads to coupled primal-dual dynamical systems. The extension of inertial acceleration to linearly constrained optimization has led to several accelerated primal–dual dynamical systems. Zeng et al. \cite{zeng_dynamical_2023} proposed the following primal–dual Nesterov-type dynamical system, which generalized the work of Su et al. \cite{su_differential_2016} to linearly constrained optimization problems. Subsequently, Bo{\c t} and Nguyen \cite{bot_improved_2021} sharpened the convergence analysis by establishing fast convergence rates for the primal–dual gap, the feasibility violation, and the objective residual. Under the additional assumption that $\nabla f$ is Lipschitz continuous, they further proved the weak convergence of the entire primal–dual trajectory. He et al. \cite{he_inertial_2023} further introduced slowly vanishing damping and general scaling coefficients into a second-order primal–dual dynamical system
\begin{equation*}
\left\{
\begin{aligned}
\ddot{x}(t)+\frac{\alpha}{t^q}\dot{x}(t)
&=
-\beta(t)\nabla_x\mathcal{L}_\rho
(
x(t),
\lambda(t)+\theta t^\kappa \dot{\lambda}(t)
)
+ e(t),
\\
\ddot{\lambda}(t)+\frac{\alpha}{t^q}\dot{\lambda}(t)
&=
\beta(t)\nabla_\lambda \mathcal{L}_\rho
(
x(t)+\theta t^\kappa \dot{x}(t),
\lambda(t)
),
\end{aligned}
\right.
\end{equation*}
where $0\le q\le\kappa\le1$, and $\beta(t)$ is a general scaling coefficient. They investigate the asymptotic properties and established fast convergence properties under different parameter regimes. Li et al. \cite{li_fast_2026} incorporated implicit Hessian damping and Tikhonov regularization into a second-order  primal–dual dynamical system and established fast convergence rates. Strong convergence of $x(t)$ to the minimum-norm solution $x^*$ is provable only when $x(t)$ eventually stays entirely inside or outside the open ball $B(0, \|x^*\|)$, if $x(t)$ repeatedly crosses the sphere, only limit inferior convergence of $x(t)$ to $x^*$ can be established. The dynamical systems above, together with several others \cite{bot_improved_2021,he_inertial_2023,zhu_fast_2025,li_fast_2026} are second-order dynamics for both the primal and dual variables. To reduce computational cost, He et al. \cite{he_second-order_2022} introduced inertial effects only in the primal variable and proposed a second-order primal + first-order dual dynamical system. They established convergence without assuming strong convexity under suitable conditions, and further proved the convergence results under different choices of the damping and scaling parameters.
Tikhonov regularization has also been incorporated into primal-dual dynamical systems. Zhu et al. \cite{zhu_fast_2025} proposed a Tikhonov-regularized inertial primal-dual dynamical system. Under appropriate conditions, Zhu et al. proved $\liminf_{t\to+\infty} \|x(t) - x^*\| = 0$. Full strong convergence of the primal trajectory was obtained only under an additional geometric condition requiring the trajectory to eventually remain either inside or outside the ball $B(0,\|\bar{x}^*\|)$\cite{zhu_fast_2025}.
Moreover, since the Tikhonov regularization acts only on the primal variable, their result concerns $x^*$ rather than the strong convergence of the primal–dual trajectory. In a different mixed-order setting, Battahi et al. \cite{fouad_battahi_simultaneous_2025} proposed a Tikhonov-regularized primal-dual dynamical system. They proved the strong convergence of the primal-dual trajectories toward a saddle point of the Lagrangian function. More recently, Zhu et al. \cite{zhu_strong_2024} studied the slowly damped inertial primal-dual dynamical system with Tikhonov regularization:
\begin{equation}
\left\{
\begin{aligned}
& \ddot{x}(t)
+\frac{\alpha}{t^q}\dot{x}(t)
+t^s
\left(
\nabla_x \mathcal{L}(x(t),\lambda(t))
+\frac{c}{t^p}x(t)
\right)
=0,
\\
& \dot{\lambda}(t)
-
t^{q+s}
\left(
\nabla_\lambda \mathcal{L}
(
x(t)+\theta t^q\dot{x}(t),
\lambda(t)
)
-\frac{c}{t^p}\lambda(t)
\right)
=0,
\end{aligned}
\right.
\end{equation}
where $0\le q<1$, $0<p<1$, $c>0$, and $s$ is a constant. Under suitable assumptions, they established strong convergence properties of the trajectories.

Despite recent advances in inertial primal-dual dynamical systems, several challenges remain. Existing Hessian-driven approaches either require explicit second-order information \cite{bot_tikhonov_2021,attouch_fast_2016,attouch_accelerated_2023} or provide only weak or conditional strong-convergence results in the settings considered \cite{alecsa_extension_2021,zhu_fast_2025}. Tikhonov regularization has been successfully used to establish strong convergence in several primal-dual dynamical systems \cite{fouad_battahi_simultaneous_2025,zhu_strong_2024,li_general_2025}. However, its combination with implicit Hessian damping in a slowly damped mixed-order primal-dual framework remains insufficiently explored. Moreover, although several dynamical systems employ power-law scaling and regularization coefficients \cite{fouad_battahi_simultaneous_2025,zhu_strong_2024}, general time-dependent coefficients have also been considered in \cite{li_general_2025}. 

Motivated by the above developments, we consider a new slowly damped inertial primal (second-order) + dual (first-order) dynamical system that incorporates implicit Hessian damping, general time scaling, and primal-dual Tikhonov regularization:
\begin{equation}\label{DS}
\begin{cases}
\ddot{x}(t) + \dfrac{\alpha}{t^q} \dot{x}(t) + \beta(t) \nabla_x \mathcal{L}_{\epsilon(t)}\bigl(x(t) +\Gamma(t)\dot{x}(t), \lambda(t)\bigr) = 0, \\
\dot{\lambda}(t) = -t^{2q} \beta(t) A\biggl[
\beta(t)\Gamma(t) \nabla_x \mathcal{L}_{\epsilon(t)}\bigl(x(t)+\Gamma(t) \dot{x}(t), \lambda(t)\bigr)
+\bigl(\frac{rq}{t^{q+1}}+\frac{\alpha-1}{t^q}\Gamma(t)\bigr)\dot{x}(t)
\biggr]\\
\qquad \quad + t^q\beta(t)\nabla_\lambda \mathcal{L}_{\epsilon(t)}\bigl(x(t) + t^q \dot{x}(t),\lambda(t)\bigr).
\end{cases}
\end{equation}
where $\Gamma(t):=r_{0}+\frac{r}{t^{q}}.$
Here \(t\geq t_{0}>0\), \(\alpha>1\), and \(0<q<1\). The parameters \(r_{0}\) and \(r\) satisfy either \(r_{0}>0\) and \(r\in\mathbb{R}\), or \(r_{0}=0\) and \(r\geq0\). Meanwhile,
\[
\mathcal{L}_{\epsilon(t)}(x,\lambda) = f(x) + \langle \lambda, Ax - b \rangle + \frac{\epsilon(t)}{2\beta(t)} \left( \|x\|^2 - \|\lambda\|^2 \right),
\]
where
\(
\beta(t) : [t_0, +\infty) \to (0, +\infty)
\)
is a non-decreasing continuously differentiable function with $\lim_{t \to +\infty} \beta(t) = +\infty$, and
\(
\epsilon(t) : [t_0, +\infty) \to (0, +\infty)
\)
is a non-increasing continuously differentiable Tikhonov regularization parameter function satisfying
\(\lim_{t \to +\infty} \epsilon(t) = 0.\) 

\subsection{Main Contributions}

\quad \ Building on the above motivation and the proposed dynamical framework, the main contributions of this work are summarized as follows.

(a) We incorporate implicit Hessian damping into a mixed-order primal–dual dynamical system with a second-order primal equation, a first-order dual equation, and slowly vanishing viscous damping ($\frac{\alpha}{t^q}$, $0 < q < 1$). Compared with existing explicit Hessian-driven damping systems \cite{attouch_first-order_2022,attouch_convergence_2023}, the proposed dynamical system effectively suppresses oscillations without requiring second-order information of $f(x)$, thereby imposing milder regularity requirements on $f$.

(b) The proposed dynamical system \eqref{DS} applies Tikhonov regularization to both the primal and dual variables and guarantees the strong convergence of the entire primal-dual trajectory to the minimum-norm primal-dual solution $\lim_{t\to+\infty}
\bigl\|(x(t),\lambda(t))-(x^{\ast},\lambda^{\ast})\bigr\|=0.$
This differs from the general \(\liminf\)-type convergence results for the primal trajectory established in \cite{zhu_fast_2025,zhu_tikhonov_2026}, for which full strong convergence requires an additional geometric condition. Full strong convergence has also been obtained in other mixed-order Tikhonov-regularized frameworks, such as \cite{fouad_battahi_simultaneous_2025,li_general_2025,zhu_strong_2024}, but those dynamical systems do not incorporate the implicit Hessian damping mechanism considered here.

(c) The time-scaling coefficient $\beta(t)$ and the Tikhonov regularization coefficient $\epsilon(t)$ in our dynamical system are both very general functions, rather than specific forms such as power functions, which broadens the applicability of the dynamical system, and we obtain fast convergence results for the primal-dual gap, the objective residual, and the feasibility violation. When $r_0 = r = 0$, $\beta(t) = t^s$, and $\epsilon(t) = c t^{s-p}$, dynamical system \eqref{DS} reduces to the dynamical system studied by Zhu et al. \cite{zhu_strong_2024} with extrapolation parameter \(\theta=1\).

\subsection{Organization}

\quad \ The remainder of this paper is organized as follows. In Section \ref{sec:prelim}, we introduce several preliminary results. In Section \ref{sec:ana}, we construct a suitable Lyapunov functional. Based on the dissipative property of the Lyapunov functional, we establish fast convergence rates for the primal-dual gap, the objective residual, and the feasibility violation. We further investigate the asymptotic behavior of the generated trajectories and prove the strong convergence of the primal-dual trajectories. In Section \ref{sec:num}, several numerical experiments are presented to illustrate the effectiveness and accelerated convergence performance of the proposed dynamical system.

\section{Preliminary results}\label{sec:prelim}

In this section, we establish auxiliary results for the convergence analysis. We first study the boundedness, strong convergence, and local
Lipschitz continuity of the regularized saddle-point trajectory. We then derive a differentiation formula for the regularized saddle-point value
and establish estimates for the feasibility violation and the Lagrangian gap.
\begin{lemma}
\label{lem2.1}
Let
\(
\Omega:=\operatorname{Saddle}(\mathcal L)
\)
be the nonempty saddle-point set of the original Lagrangian
$\mathcal L$, and let
\(
y^*:=(x^*,\lambda^*)=P_{\Omega}(0_{\mathcal X},0_{\mathcal Y})
\)
be the minimum-norm element of \ $\Omega$.
Assume that
\[
\frac{\epsilon(t)}{\beta(t)}\in C^1([t_0,+\infty);(0,+\infty))
\ \text{and}\ \lim_{t\to+\infty}\frac{\epsilon(t)}{\beta(t)}=0.
\]
For every $t\ge t_0$, let 
\( 
(x_t,\lambda_t) 
\)
be the unique saddle point of $\mathcal{L}_{\epsilon(t)}$. Then
\begin{enumerate}
\item[(i)] for all $t>t_{0}$, $\norm{(x_t,\lambda_t)}\leq\norm{(x^{*}, \lambda^{*})}$ and $\lim _{t \to +\infty}\norm{(x_t, \lambda_t)-(x^{*}, \lambda^{*})}=0$,
    \item[(ii)] the mapping
\(
t\mapsto (x_t,\lambda_t)
\)
is locally Lipschitz continuous on $(t_0,+\infty)$ and hence locally
absolutely continuous. Moreover, for a.e. $t>t_0$,
\[
\|(\dot x_t,\dot\lambda_t)\| \le \left|\left(\frac{\epsilon(t)}{\beta(t)}\right)' \right| \frac{\beta(t)}{\epsilon(t)}\;\|(x_t,\lambda_t)\|
     \le\left|\left(\frac{\epsilon(t)}{\beta(t)}\right)' \right|  \frac{\beta(t)}{\epsilon(t)}\;\|(x^*,\lambda^*)\|.
\]
\end{enumerate}
\end{lemma}

\begin{proof}
(i) Consider the operator $\mathscr{M}: \mathcal{X} \times \mathcal{Y} \to \mathcal{X} \times \mathcal{Y}$ defined by
\[
{\mathscr{M}(x, \lambda):=(\nabla_{x} \mathcal{L}(x, \lambda),-\nabla_{\lambda} \mathcal{L}(x, \lambda))}.
\]
Since $\mathcal{L}$ is convex in $x$ and concave in $\lambda$, $M$ is monotone. Moreover, $M$ is continuous on $\mathcal{X}\times\mathcal{Y}$, and hence maximally monotone. Its zero set coincides with the saddle-point set of $\mathcal{L}$, namely,
\(
\operatorname{zer}M=\Omega.
\)
Let $w(t):=(x_t,\lambda_t)$. From the optimality conditions for the saddle point of $\mathcal{L}_{\epsilon(t)}$,
\[
M(w(t))+\frac{\epsilon(t)}{\beta(t)}w(t)=0.
\]
Since $y^*:=(x^*,\lambda^*)=P_{\Omega}(0_{\mathcal{X}},0_{\mathcal{Y}})$ satisfies $M(y^*)=0$, the monotonicity of $M$ gives
\[
0\le \bigl\langle M(w(t))-M(y^*),w(t)-y^*\bigr\rangle =-\frac{\epsilon(t)}{\beta(t)}\bigl\langle w(t),w(t)-y^*\bigr\rangle.
\]
Hence,
$\|w(t)\|^2\le \bigl\langle w(t),y^*\bigr\rangle\le\|w(t)\|\,\|y^*\|,$
and therefore
$\|w(t)\|\le\|y^*\|.$
Furthermore,
\[
0\in M(w(t))+\frac{\epsilon(t)}{\beta(t)}\big(w(t)-(0_{\mathcal{X}},0_{\mathcal{Y}})\big).
\]
Since
$\frac{\epsilon(t)}{\beta(t)}\to 0 \ \text{as } t\to+\infty,$
for sufficiently large $t$, $\frac{\epsilon(t)}{\beta(t)}\in(0,1)$, so \cite{bauschke2017convex} 
Theorem 23.44 is applied to $M$ and $x=(0_{\mathcal{X}},0_{\mathcal{Y}})$ yields
\[
w(t)\to P_{\operatorname{zer}M}(0_{\mathcal{X}},0_{\mathcal{Y}})=P_{\Omega}(0_{\mathcal{X}},0_{\mathcal{Y}})=y^*
\]
strongly as $t\to+\infty$. Thus,
\[
\lim_{t\to+\infty}\bigl\|(x_t,\lambda_t)-(x^*,\lambda^*)\bigr\|=0.
\]

(ii) From \eqref{gradient of L 1} and \eqref{gradient of L 2}, for $t>t_{0}$ and $h$ near zero
\[
\mathscr{M}(w(t)) = -\frac{\epsilon(t)}{\beta(t)}w(t),\qquad \mathscr{M}(w(t+h)) = -\frac{\epsilon(t+h)}{\beta(t+h)}w(t+h).
\]
Monotonicity of $\mathscr{M}$ gives
\[
\langle \mathscr{M}(w(t+h))-\mathscr{M}(w(t)),\; w(t+h)-w(t)\rangle \ge 0.
\]
Substituting the expressions for $\mathscr{M}$ yields
\[
\bigl\langle -\frac{\epsilon(t+h)}{\beta(t+h)}w(t+h)+\frac{\epsilon(t)}{\beta(t)}w(t),\; w(t+h)-w(t)\bigr\rangle \ge 0.
\]
Rearranging,
\[
{\frac{\epsilon(t+h)}{\beta(t+h)}\,\|w(t+h)-w(t)\|^2 \le \bigl(\frac{\epsilon(t)}{\beta(t)}-\frac{\epsilon(t+h)}{\beta(t+h)}\bigr)\langle w(t),\, w(t+h)-w(t)\rangle}.
\]
Take absolute values and apply the Cauchy–Schwarz inequality:
\[
\frac{\epsilon(t+h)}{\beta(t+h)}\,\|w(t+h)-w(t)\|^2 \le \left|\frac{\epsilon(t+h)}{\beta(t+h)}-\frac{\epsilon(t)}{\beta(t)}\right|\;\|w(t)\|\;\|w(t+h)-w(t)\|.
\]
If $\|w(t+h)-w(t)\|>0$ we can cancel one factor to obtain
\begin{equation}\label{inner}
\|w(t+h)-w(t)\| \le \frac{\left|\frac{\epsilon(t+h)}{\beta(t+h)}-\frac{\epsilon(t)}{\beta(t)}\right|}{\frac{\epsilon(t+h)}{\beta(t+h)}}\;\|w(t)\|.     
\end{equation}
Now use the differentiability of $\frac{\epsilon(t)}{\beta(t)}$. By the mean value theorem, for some $\xi\in(0,1)$,
\[
\frac{\epsilon(t+h)}{\beta(t+h)}-\frac{\epsilon(t)}{\beta(t)} = \left(\frac{\epsilon(t+\xi h)}{\beta(t+\xi h)}\right)'h.
\]
Thus
\[
\frac{\left|\dfrac{\epsilon(t+h)}{\beta(t+h)}-\dfrac{\epsilon(t)}{\beta(t)}\right|}{\dfrac{\epsilon(t+h)}{\beta(t+h)}}
=
\frac{\left| \left( \dfrac{\epsilon(t+\xi h)}{\beta(t+\xi h)} \right)' \right|}{\dfrac{\epsilon(t+\xi h)}{\beta(t+\xi h)}}
\cdot
\frac{\dfrac{\epsilon(t+\xi h)}{\beta(t+\xi h)}}{\dfrac{\epsilon(t+h)}{\beta(t+h)}}
\, |h|.
\]
Considering separately $h \to 0^+$ and $h \to 0^-$, and using the continuity of $t \mapsto (x_t, \lambda_t)$, we obtain $\frac{\epsilon(t+\xi h)}{\beta(t+\xi h)} {\frac{\beta(t+h)}{\epsilon(t+h)}}\to1$ and $\left|\left(\frac{\epsilon(t+\xi h)}{\beta(t+\xi h)}\right)'\right|  \frac{\beta(t+\xi h)}{\epsilon(t+\xi h)}\to\left|\left(\frac{\epsilon(t)}{\beta(t)}\right)' \right| \frac{\beta(t)}{\epsilon(t)}$ . Consequently
\[
\frac{\left|\frac{\epsilon(t+h)}{\beta(t+h)}-\frac{\epsilon(t)}{\beta(t)}\right|}{\frac{\epsilon(t+h)}{\beta(t+h)}} = \left|\left(\frac{\epsilon(t)}{\beta(t)}\right)'\right| \frac{\beta(t)}{\epsilon(t)}\,|h| + o(|h|).
\]
We get the trajectory $w(t)$ is Lipschitz continuous on each compact interval in $(t_0, +\infty)$. 
Consequently, $w$ is locally Lipschitz continuous, and therefore
locally absolutely continuous on $(t_0,+\infty)$. In particular,
$w$ is differentiable for a.e. $t>t_0$.
Insert this into \eqref{inner}, divide by $|h|$ and let $h\to0$. The left‑hand side becomes  (for almost every $t$), and we get
\[
\|\dot w(t)\| \le \left|\left(\frac{\epsilon(t)}{\beta(t)}\right)'\right| \frac{\beta(t)}{\epsilon(t)}\;\|w(t)\| \quad\text{a.e.}
\]
Using (i), $\|w(t)\|\le\|(x^*,\lambda^*)\|$,  meaning that (ii) is satisfied.
\end{proof}

\begin{lemma}
Let $(x_t,\lambda_t)$ denote the unique saddle point of
$\mathcal{L}_{\epsilon(t)}$. Then for each $t>t_0$, 
\[
\frac{d}{dt}\mathcal{L}_{\epsilon(t)}(x_t,\lambda_t)
=
\left(\frac{\epsilon(t)}{2\beta(t)}\right)'
\left(\|x_t\|^2-\|\lambda_t\|^2\right).
\]
\end{lemma}

\begin{proof}
Fix $t>t_0$ and let $h$ be sufficiently small. Since $(x_t,\lambda_t)$ is the saddle point of $\mathcal{L}_{\epsilon(t)}$, 
\[
\mathcal{L}_{\epsilon(t)}(x_t,\lambda_{t+h})
\le
\mathcal{L}_{\epsilon(t)}(x_t,\lambda_t)
\le
\mathcal{L}_{\epsilon(t)}(x_{t+h},\lambda_t).
\]
On the other hand, since $(x_{t+h},\lambda_{t+h})$ is the saddle point of $\mathcal{L}_{\epsilon(t+h)}$, we also have
\[
-
\mathcal{L}_{\epsilon(t+h)}(x_t,\lambda_{t+h})
\le
-
\mathcal{L}_{\epsilon(t+h)}(x_{t+h},\lambda_{t+h})
\le
-
\mathcal{L}_{\epsilon(t+h)}(x_{t+h},\lambda_t).
\]
Adding the two inequalities yields
\[
\mathcal{L}_{\epsilon(t)}(x_t,\lambda_t)
-
\mathcal{L}_{\epsilon(t+h)}(x_{t+h},\lambda_{t+h})
\le
\mathcal{L}_{\epsilon(t)}(x_{t+h},\lambda_t)
-
\mathcal{L}_{\epsilon(t+h)}(x_{t+h},\lambda_t).
\]
Using the definition of $\mathcal{L}_{\epsilon(t)}$, we compute
\[
\begin{aligned}
&
\mathcal{L}_{\epsilon(t)}(x_{t+h},\lambda_t)
-
\mathcal{L}_{\epsilon(t+h)}(x_{t+h},\lambda_t)
\\
&=
\frac{\epsilon(t)}{2\beta(t)}
\left(\|x_{t+h}\|^2-\|\lambda_t\|^2\right)
-
\frac{\epsilon(t+h)}{2\beta(t+h)}
\left(\|x_{t+h}\|^2-\|\lambda_t\|^2\right)
\\
&=
\left(
\frac{\epsilon(t)}{2\beta(t)}
-
\frac{\epsilon(t+h)}{2\beta(t+h)}
\right)
\left(\|x_{t+h}\|^2-\|\lambda_t\|^2\right).
\end{aligned}
\]
Since $\epsilon(t)$ and $\beta(t)$ are differentiable, so
\[
\frac{\epsilon(t+h)}{2\beta(t+h)}
=
\frac{\epsilon(t)}{2\beta(t)}
+
h
\left(\frac{\epsilon(t)}{2\beta(t)}\right)'
+
o(h).
\]
Therefore,
\[
\mathcal{L}_{\epsilon(t)}(x_{t+h},\lambda_t)
-
\mathcal{L}_{\epsilon(t+h)}(x_{t+h},\lambda_t)
=
-
h
\left(\frac{\epsilon(t)}{2\beta(t)}\right)'
\left(\|x_{t+h}\|^2-\|\lambda_t\|^2\right)
+
o(h).
\]
Similarly, we obtain
\[
\mathcal{L}_{\epsilon(t)}(x_t,\lambda_t)
-
\mathcal{L}_{\epsilon(t+h)}(x_{t+h},\lambda_{t+h})
\ge
-
h
\left(\frac{\epsilon(t)}{2\beta(t)}\right)'
\left(\|x_t\|^2-\|\lambda_{t+h}\|^2\right)
+
o(h).
\]
By Lemma \ref{lem2.1} (ii), the mapping $t\mapsto (x_t,\lambda_t)$ is locally Lipschitz continuous. Hence,
\[
x_{t+h}\to x_t
\ \text{and}\
\lambda_{t+h}\to\lambda_t
\ \text{as }h\to0.
\]
Dividing both sides by $h$ and letting $h\to0$, we obtain
\[
\frac{d}{dt}\mathcal{L}_{\epsilon(t)}(x_t,\lambda_t)
=
\left(\frac{\epsilon(t)}{2\beta(t)}\right)'
\left(\|x_t\|^2-\|\lambda_t\|^2\right).
\]
This completes the proof.
\end{proof}

\begin{lemma} \label{lem2.3}
Let $y^*=(x^*,\lambda^*)$, $y_t=(x_t,\lambda_t)$ and $(x^*,\lambda^*)= P_{\Omega}(0_{\mathcal X},0_{\mathcal Y})$. Let $(x,\lambda)$ be  the trajectory generated by the dynamical system \eqref{DS} and let $(x_t,\lambda_t)$ be the unique saddle point of $\mathcal{L}_{\epsilon(t)}(x,\lambda)$. Then, for every $t \ge t_0$, the following two estimates are satisfied:
\begin{enumerate}
  \item[(i)] 
  \begin{equation}
    \|Ax - b\| \leq \|A\| \ \|x - x_t\| + \frac{\epsilon(t)}{\beta(t)} \|(x^*,\lambda ^*)\|,
  \end{equation}
\item[(ii)]
\begin{equation}
\begin{split}
\mathcal{L}\bigl(x(t) + \Gamma(t)\dot{x}(t), \lambda^*\bigr) - \mathcal{L}(x^*, \lambda^*)
\le \mathcal{L}_{\epsilon(t)}\bigl(x(t) + \Gamma(t)\dot{x}(t), \lambda_t\bigr) - \mathcal{L}_{\epsilon(t)}(x_t, \lambda_t) \\
\quad + \|\lambda_t - \lambda^*\| \left\| A\bigl(x(t) + \Gamma(t)\dot{x}(t)\bigr) - b \right\|
+ \frac{\epsilon(t)}{2\beta(t)} \left( \|x^*\|^2 - \bigl\|x(t) + \Gamma(t)\dot{x}(t)\bigr\|^2 \right).
\end{split}
\end{equation}
\end{enumerate}
\end{lemma}

\begin{proof}
(i) By the definition of the regularized Lagrangian, it follows that
\begin{equation}\label{gradient of L 1}
\begin{split}
\nabla_x \mathcal{L}_{\epsilon(t)}(x, \lambda) &= \nabla f(x) + A^\top \lambda + \dfrac{\epsilon(t)}{\beta(t)} x 
\end{split}
\end{equation}
\begin{equation}\label{gradient of L 2}
\begin{split}
\nabla_\lambda \mathcal{L}_{\epsilon(t)}(x, \lambda) &= A x - b - \dfrac{\epsilon(t)}{\beta(t)} \lambda,
\end{split}
\end{equation}
the saddle point $(x_t, \lambda_t)$ satisfies the dynamical system
\[
\begin{cases}
\nabla_x \mathcal{L}_{\epsilon(t)}(x_t, \lambda_t) = \nabla f(x_t) + A^\top \lambda_t + \dfrac{\epsilon(t)}{\beta(t)} x_t=0, \\[1ex]
\nabla_\lambda \mathcal{L}_{\epsilon(t)}(x_t, \lambda_t) = Ax_t - b - \dfrac{\epsilon(t)}{\beta(t)} \lambda_t=0.
\end{cases}
\]
With the conclusion in Lemma \ref{lem2.1}, we provide the following norm estimate:
\begin{equation*}
\begin{split}
\|Ax - b\| &= \|A(x - x_t) + A x_t - b\| \\
&\le \|A\| \ \|x - x_t\| + \|A x_t - b\| \\
&= \|A\| \ \|x - x_t\| + \frac{\epsilon(t)}{\beta(t)} \|\lambda_t\| \\
&\le \|A\| \ \|x - x_t\| + \frac{\epsilon(t)}{\beta(t)} \|y^*\|,
\end{split}
\end{equation*}

(ii) In addition, from the property $\mathcal{L}_{\epsilon(t)}(x_t, \lambda_t) \leq \mathcal{L}_{\epsilon(t)}(x^*, \lambda_t)$ and the fact that $(x^*, \lambda^*) \in \Omega$,we obtain
\begin{align*}
& \mathcal{L}_{\epsilon(t)}\big(x(t) + \Gamma(t)\dot{x}(t), \lambda_t\big) - \mathcal{L}_{\epsilon(t)}(x_t, \lambda_t) \\
& \quad \geq \mathcal{L}_{\epsilon(t)}\big(x(t) + \Gamma(t)\dot{x}(t), \lambda_t\big) - \mathcal{L}_{\epsilon(t)}(x^*, \lambda_t) \\
& \quad = \mathcal{L}\big(x(t) + \Gamma(t)\dot{x}(t), \lambda_t\big) - \mathcal{L}(x^*, \lambda_t) + \frac{\epsilon(t)}{2\beta(t)} \left(\|x(t)+ \Gamma(t)\dot{x}(t)\|^2 - \|x^*\|^2\right) \\
& \quad = \mathcal{L}\big(x(t) + \Gamma(t)\dot{x}(t), \lambda_t\big) - \mathcal{L}(x^*, \lambda^*) + \frac{\epsilon(t)}{2\beta(t)} \left(\|x(t)+ \Gamma(t)\dot{x}(t)\|^2 - \|x^*\|^2\right) \\
& \quad = \mathcal{L}\big(x(t) + \Gamma(t)\dot{x}(t), \lambda^*\big) - \mathcal{L}(x^*, \lambda^*) \\
& \qquad + \left\langle \lambda_t - \lambda^*, A\big(x(t) + \Gamma(t)\dot{x}(t)\big) - b \right\rangle + \frac{\epsilon(t)}{2\beta(t)} \left(\|x(t)+ \Gamma(t)\dot{x}(t)\|^2 - \|x^*\|^2\right) \\
& \quad \geq \mathcal{L}\big(x(t) + \Gamma(t)\dot{x}(t), \lambda^*\big) - \mathcal{L}(x^*, \lambda^*) \\
& \qquad - \|\lambda_t - \lambda^*\| \left\| A\big(x(t) + \Gamma(t)\dot{x}(t)\big) - b \right\| + {\frac{\epsilon(t)}{2\beta(t)} \left(\|x(t)+ \Gamma(t)\dot{x}(t)\|^2-\|x^*\|^2 \right)},
\end{align*}
Rearranging the above inequality yields
\begin{equation*}
\begin{split}
\mathcal{L}\bigl(x(t) + \Gamma(t)\dot{x}(t), \lambda^*\bigr) - \mathcal{L}(x^*, \lambda^*)
\le \mathcal{L}_{\epsilon(t)}\bigl(x(t) + \Gamma(t)\dot{x}(t), \lambda_t\bigr) - \mathcal{L}_{\epsilon(t)}(x_t, \lambda_t) \\
\quad + \|\lambda_t - \lambda^*\| \left\| A\bigl(x(t) + \Gamma(t)\dot{x}(t)\bigr) - b \right\|
+ \frac{\epsilon(t)}{2\beta(t)} \left( \|x^*\|^2 - \bigl\|x(t) + \Gamma(t)\dot{x}(t)\bigr\|^2 \right).
\end{split}
\end{equation*}
This completes the proof.
\end{proof}

\section{Convergence Analysis}\label{sec:ana}

In this section, we first establish a differential inequality for a Lyapunov functional. Under additional assumptions on the regularization parameter, we then derive convergence rates and prove strong convergence of the primal-dual trajectory to the minimum-norm saddle point.

\begin{theorem}\label{Th3.1}
Let $(x,\lambda): [t_0, +\infty) \to \mathcal{X}\times\mathcal{Y}$ be any solution of \eqref{DS} and $\alpha>1$. Define
\begin{equation}
\begin{split}
\mathcal{E}(t) &= t^{2q} \beta(t) \bigl( \mathcal{L}_{\epsilon(t)}\bigl(
x(t) + \bigl(r_0 + \dfrac{r}{t^q}\bigr) \dot{x}(t), \lambda_t
\bigr) - \mathcal{L}_{\epsilon(t)}\bigl(x_t, \lambda_t\bigr) \bigr) + \dfrac{1}{2} \bigl\| x(t) - x_t + t^q \dot{x}(t) \bigr\|^2 \\
&\quad + \dfrac{1}{2} \bigl\| \lambda(t) - \lambda_t \bigr\|^2 + \dfrac{\alpha - 1 - q t^{q-1}}{2} \bigl\| x(t) - x_t \bigr\|^2.
\label{eq:energy}
\end{split}
\end{equation}
Let the constant W satisfy $0 < W < \min\left\{ \frac{\alpha}{2},\ \frac{\alpha}{2(\alpha+1)} \right\}$, assume that the following conditions $H_0$ hold:
\[
\begin{cases} \tag{$H_0$}\label{H0}
\displaystyle\lim_{t\to+\infty} \frac{t^{q-1}\beta(t)}{\epsilon(t)} = 0,\\
\epsilon(t)>\dfrac{\alpha}{t^{q+m}} \quad \text{for all sufficiently large } t,\\
\displaystyle\limsup_{t\to+\infty}t^q\frac{\dot\beta(t)}{\beta(t)}<1,\\
0< q < m <1.
\end{cases}
\]
then there exists $T_0 \ge t_0$ and constants $P_i \ge 0, i=1,...,5$, such that for a.e. $t \ge T$
\begin{equation} \begin{split}\label{Compute}
\dot{\mathcal{E}}(t) + \frac{W}{t^m} \mathcal{E}(t)
\le{} & - P_1 t^{q} \beta(t) \bigl(\mathcal{L}_{\epsilon(t)}(x(t)+\Gamma(t) \dot{x}(t), \lambda_t) - \mathcal{L}_{\epsilon(t)}(x_t, \lambda_t)\bigr)
- P_2 \epsilon(t) t^q \|x(t)-x_t\|^2\\ 
&-P_3 t^q \|\dot{x}(t)\|^2 - P_4 \epsilon(t) t^q \|\lambda(t) - \lambda_t\|^2 \\
&- P_5 t^{q} \beta(t)^2 \|\nabla_x \mathcal{L}_{\epsilon(t)}(x(t)+\Gamma(t) \dot{x}(t), \lambda(t))\|^2 + \left( \frac{2\alpha^2}{ t^q \epsilon(t)} + \frac{t^q}{\alpha-1} \right) \|\dot{x}_t\|^2\\
&+ \left( \frac{r_0+1}{2} t^{2q+1} \beta(t) + \frac{r}{2} t^{q+1} \beta(t) + \frac{1}{t^q \epsilon(t)} \right) \|\dot{\lambda}_t\|^2\\
&+ t^{2q} \beta(t)\left(\frac{\epsilon(t)}{2\beta(t)}\right)' (\|x(t)+\Gamma(t) \dot{x}(t)\|^2-\|x_t\|^2). 
\end{split} \end{equation}
\end{theorem}

\begin{proof}

We first differentiate each component of the Lyapunov functional separately, writing ${\mathcal{E}}(t)$ = ${\mathcal{E}_1}(t)$ + ${\mathcal{E}_2}(t)$ + ${\mathcal{E}_3}(t)$ + ${\mathcal{E}_4}(t)$,
where
\[
{\mathcal{E}_1}(t) = t^{2q} \beta(t) \left( \mathcal{L}_{\epsilon(t)}\bigl( x(t) + \Gamma(t) \dot{x}(t), \lambda_t \bigr) - \mathcal{L}_{\epsilon(t)}(x_t, \lambda_t) \right),
\]
\[
{\mathcal{E}_2}(t) = \frac{1}{2} \left\| x(t) - x_t + t^q \dot{x}(t) \right\|^2,
\]
\[
{\mathcal{E}_3}(t) = \frac{1}{2} \left\| \lambda(t) - \lambda_t \right\|^2,
\]
and
\[
{\mathcal{E}_4}(t) = \frac{\alpha - 1 - q t^{q-1}}{2} \left\| x(t) - x_t \right\|^2.
\]
Differentiating the Lyapunov functional ${\mathcal{E}}(t)$, this yields
\begin{equation} \begin{split}
\dot{\mathcal{E}}_1(t) &= \bigl(2q t^{2q-1} \beta(t) + t^{2q} \dot{\beta}(t)\bigr)\bigl(\mathcal{L}_{\epsilon(t)}(x(t)+\Gamma(t) \dot{x}(t), \lambda_t) - \mathcal{L}_{\epsilon(t)}(x_t, \lambda_t)\bigr) \\
&\quad + t^{2q} \beta(t) \ \Bigl[ \langle \nabla_x \mathcal{L}_{\epsilon(t)}(x(t)+\Gamma(t) \dot{x}(t), \lambda_t), \dot{x}(t) + \dot{\Gamma}(t) \dot{x}(t) + \Gamma(t) \ddot{x}(t) \rangle \\
&\quad + \langle \nabla_\lambda \mathcal{L}_{\epsilon(t)}(x(t)+\Gamma(t) \dot{x}(t), \lambda_t), \dot{\lambda}_t \rangle - \langle \nabla_x \mathcal{L}_{\epsilon(t)}(x_t, \lambda_t), \dot{x}_t \rangle - \langle \nabla_\lambda \mathcal{L}_{\epsilon(t)}(x_t, \lambda_t), \dot{\lambda}_t \rangle \\
&\quad + \left(\frac{\epsilon(t)}{2\beta(t)}\right)' \bigl(\|x(t)+\Gamma(t) \dot{x}(t)\|^2 - \|\lambda_t\|^2\bigr) - \left(\frac{\epsilon(t)}{2\beta(t)}\right)' \bigl(\|x_t\|^2 - \|\lambda_t\|^2\bigr) \Bigr].
\end{split} \end{equation}
Using the gradient of Tikhonov-regularized Lagrangian function $\nabla_\lambda \mathcal{L}_{\epsilon(t)}(x_t, \lambda_t) =0$
\begin{equation*} \begin{split}
\langle \nabla_\lambda \mathcal{L}_{\epsilon(t)}(x(t)+\Gamma(t) \dot{x}(t), \lambda_t), \dot{\lambda}_t \rangle &= \left\langle A(x(t)+\Gamma(t) \dot{x}(t)) - b - \frac{\epsilon(t)}{\beta(t)} \lambda_t, \dot{\lambda}_t \right\rangle \\
&= \left\langle A(x(t)+\Gamma(t) \dot{x}(t) - x_t) , \dot{\lambda}_t \right\rangle,
\end{split} \end{equation*}
we obtain
\begin{align*}
\dot{\mathcal{E}}_1
&= \bigl(2q t^{2q-1} \beta(t) + t^{2q} \dot{\beta}(t)\bigr)\bigl(\mathcal{L}_{\epsilon(t)}(x(t)+\Gamma(t) \dot{x}(t), \lambda_t) - \mathcal{L}_{\epsilon(t)}(x_t, \lambda_t)\bigr) \\
&\quad + t^{2q} \beta(t)  \langle \nabla_x \mathcal{L}_{\epsilon(t)}(x(t)+\Gamma(t) \dot{x}(t), \lambda_t), \dot{x}(t) + \dot{\Gamma}(t) \dot{x}(t) + \Gamma(t) \ddot{x}(t) \rangle \\
&\quad + t^{2q} \beta(t) \left\langle A(x(t)+\Gamma(t) \dot{x}(t) - x_t) , \dot{\lambda}_t \right\rangle \\
&\quad + t^{2q} \beta(t)\left(\frac{\epsilon(t)}{2\beta(t)}\right)' \bigl(\|x(t)+\Gamma(t) \dot{x}(t)\|^2 - \|x_t\|^2\bigr) \\
&= \bigl(2q t^{2q-1} \beta(t) + t^{2q} \dot{\beta}(t)\bigr)\bigl(\mathcal{L}_{\epsilon}(x(t)+\Gamma(t) \dot{x}(t),\lambda_t) - \mathcal{L}_{\epsilon}(x_t,\lambda_t)\bigr) \\
&\quad + t^{2q} \beta(t) (1+\dot{\Gamma}(t)) \ \langle \nabla_x \mathcal{L}_{\epsilon}(x(t)+\Gamma(t) \dot{x}(t),\lambda_t), \dot{x}(t) \rangle \\
&\quad - \alpha t^q \beta(t) \Gamma(t) \langle \nabla_x \mathcal{L}_{\epsilon}(x(t)+\Gamma(t) \dot{x}(t),\lambda_t), \dot{x}(t) \rangle \\
&\quad - t^{2q} \beta(t)^2 \Gamma(t) \langle \nabla_x \mathcal{L}_{\epsilon}(x(t)+\Gamma(t) \dot{x}(t),\lambda_t), \nabla_x \mathcal{L}_{\epsilon}(x(t)+\Gamma(t) \dot{x}(t),\lambda(t)) \rangle \\
&\quad + t^{2q} \beta(t) \ \langle A(x(t)+\Gamma(t) \dot{x}(t)-x_t) , \dot{\lambda}_t \rangle \\
&\quad + t^{2q} \beta(t) \left(\frac{\epsilon(t)}{2\beta(t)}\right)' \bigl(\|x(t)+\Gamma(t) \dot{x}(t)\|^2 - \|x_t\|^2\bigr),
\end{align*}
where we use the equation of motion in \eqref{DS} to substitute the acceleration term $\ddot{x}(t)$.\\
Using the affine dependence on the dual variable, this yields
\begin{equation*}
\begin{aligned}
&- t^{2q} \beta(t)^2 \Gamma(t) \langle \nabla_x \mathcal{L}_{\epsilon(t)}(x(t)+\Gamma(t) \dot{x}(t), \lambda_t), \nabla_x \mathcal{L}_{\epsilon(t)}(x(t)+\Gamma(t) \dot{x}(t), \lambda(t)) \rangle \\
&= - t^{2q} \beta(t)^2 \Gamma(t) \langle \nabla_x \mathcal{L}_{\epsilon(t)}(x(t)+\Gamma(t) \dot{x}(t), \lambda(t)) - A^\top(\lambda(t) - \lambda_t), \nabla_x \mathcal{L}_{\epsilon(t)}(x(t)+\Gamma(t) \dot{x}(t), \lambda(t)) \rangle \\
&= - t^{2q} \beta(t)^2 \Gamma(t) \|\nabla_x \mathcal{L}_{\epsilon(t)}(x(t)+\Gamma(t) \dot{x}(t), \lambda(t))\|^2 \\
&\quad + t^{2q} \beta(t)^2 \Gamma(t) \langle A \nabla_x \mathcal{L}_{\epsilon(t)}(x(t)+\Gamma(t) \dot{x}(t), \lambda(t)), (\lambda(t) - \lambda_t) \rangle,
\end{aligned}
\end{equation*}
and 
\begin{equation*}
\begin{split}
t^{2q} \beta(t)  \left\langle A(x(t)+\Gamma(t) \dot{x}(t) - x_t) , \dot{\lambda}_t \right\rangle = t^{2q} \beta(t)  \langle A(x(t) - x_t), \dot{\lambda}_t \rangle + t^{2q} \beta(t) \Gamma(t) \langle A \dot{x}(t), \dot{\lambda}_t \rangle.
\end{split}
\end{equation*}
Let \( u(t) = x(t) - x_t + t^q \dot{x}(t) \), then
\begin{equation*} \begin{split}
\dot{u}(t) &= \dot{x}(t) - \dot{x}_t + q t^{q-1} \dot{x}(t) + t^q \ddot{x}(t) \\
&= \bigl(1 + q t^{q-1}\bigr) \dot{x}(t) - \dot{x}_t - t^q \left[ \frac{\alpha}{t^q} \dot{x}(t) + \beta(t) \nabla_x \mathcal{L}_{\epsilon(t)}(x(t)+\Gamma(t) \dot{x}(t), \lambda(t)) \right] \\
&= \bigl(1 + q t^{q-1} - \alpha\bigr) \dot{x}(t) - \dot{x}_t - t^q \beta(t) \nabla_x \mathcal{L}_{\epsilon(t)}(x(t)+\Gamma(t) \dot{x}(t), \lambda_t) - t^q \beta(t) A^\top(\lambda(t) - \lambda_t).
\end{split} \end{equation*}
so
\begin{equation} \label{dE2}\begin{split}
\dot{\mathcal{E}}_2(t) &= \langle u(t), \dot{u}(t) \rangle \\
&= \langle x(t) - x_t + t^q \dot{x}(t), \bigl(1 + q t^{q-1} - \alpha\bigr) \dot{x}(t) - \dot{x}_t - t^q \beta(t) \nabla_x \mathcal{L}_{\epsilon(t)}(x(t)+\Gamma(t) \dot{x}(t), \lambda_t) \\
&\quad - t^q \beta(t) A^\top(\lambda(t) - \lambda_t) \rangle \\
&= \bigl(1 + q t^{q-1} - \alpha\bigr) \langle x(t) - x_t, \dot{x}(t) \rangle + t^q \bigl(1 + q t^{q-1} - \alpha\bigr) \|\dot{x}(t)\|^2\ - \langle x(t) - x_t, \dot{x}_t \rangle  \\
&\quad- t^q \beta(t) \langle x(t) - x_t, \nabla_x \mathcal{L}_{\epsilon(t)}(x(t)+\Gamma(t) \dot{x}(t), \lambda_t) \rangle\\
&\quad - t^{2q} \beta(t) \langle \dot{x}(t), \nabla_x \mathcal{L}_{\epsilon(t)}(x(t)+\Gamma(t) \dot{x}(t), \lambda_t) \rangle \\
&\quad- t^{q} \beta(t) \langle x(t) - x_t, A^\top(\lambda(t) - \lambda_t) \rangle - t^{2q} \beta(t) \langle \dot{x}(t), A^\top(\lambda(t) - \lambda_t) \rangle- t^q \langle \dot{x}(t), \dot{x}_t \rangle.
\end{split} \end{equation}
Since $ \mathcal{L}_{\epsilon(t)}(\cdot, \lambda_t) $ is $\frac{\epsilon(t)}{\beta(t)}$-strongly convex with respect to the primal variable, we can invoke the standard strong convexity inequality to estimate the term 
\begin{equation*} \begin{split}
&- \langle x(t) - x_t, \nabla_x \mathcal{L}_{\epsilon(t)}(x(t)+\Gamma(t) \dot{x}(t), \lambda_t) \rangle \\
&=- \langle x(t)+\Gamma(t)\dot{x}(t)- x_t, \nabla_x \mathcal{L}_{\epsilon(t)}(x(t)+\Gamma(t) \dot{x}(t), \lambda_t) \rangle +\langle \Gamma(t)\dot{x}(t), \nabla_x \mathcal{L}_{\epsilon(t)}(x(t)+\Gamma(t) \dot{x}(t), \lambda_t) \rangle \\
&\leq - \Bigl[ \mathcal{L}_{\epsilon(t)}(x(t)+\Gamma(t) \dot{x}(t), \lambda_t) - \mathcal{L}_{\epsilon(t)}(x_t, \lambda_t)\Bigr] - \frac{\epsilon(t)}{2\beta(t)} \|x(t)+\Gamma(t) \dot{x}(t) - x_t\|^2 \\
&\quad +\langle \Gamma(t)\dot{x}(t), \nabla_x \mathcal{L}_{\epsilon(t)}(x(t)+\Gamma(t) \dot{x}(t), \lambda_t) \rangle  \\
&= - \Bigl[ \mathcal{L}_{\epsilon(t)}(x(t)+\Gamma(t) \dot{x}(t), \lambda_t) - \mathcal{L}_{\epsilon(t)}(x_t, \lambda_t) \Bigr] - \frac{\epsilon(t)}{2\beta(t)} \|x(t) - x_t\|^2 - \frac{\epsilon(t)}{\beta(t)} \langle x(t) - x_t, \Gamma(t) \dot{x}(t) \rangle \\
&\quad - \frac{\epsilon(t)}{2\beta(t)} \|\Gamma(t) \dot{x}(t)\|^2+\langle \Gamma(t)\dot{x}(t), \nabla_x \mathcal{L}_{\epsilon(t)}(x(t)+\Gamma(t) \dot{x}(t), \lambda_t) \rangle .
\end{split} \end{equation*}
This yields
\begin{equation*} \begin{split}
&- t^q \beta(t) \langle x(t) - x_t, \nabla_x \mathcal{L}_{\epsilon(t)}(x(t)+\Gamma(t) \dot{x}(t), \lambda_t) \rangle \\
&\leq - t^q \beta(t) \Bigl[ \mathcal{L}_{\epsilon(t)}(x(t)+\Gamma(t) \dot{x}(t), \lambda_t) - \mathcal{L}_{\epsilon(t)}(x_t, \lambda_t) \Bigr] - t^q \beta(t) \frac{\epsilon(t)}{2\beta(t)} \|x(t) - x_t\|^2 \\
&\quad - t^q \beta(t) \frac{\epsilon(t)}{\beta(t)} \Gamma(t) \langle x(t) - x_t, \dot{x}(t) \rangle - t^q \beta(t) \frac{\epsilon(t)}{2\beta(t)} \Gamma(t)^2 \|\dot{x}(t)\|^2\\
&\quad +t^q\beta(t)\langle \Gamma(t)\dot{x}(t), \nabla_x \mathcal{L}_{\epsilon(t)}(x(t)+\Gamma(t) \dot{x}(t), \lambda_t) \rangle \\
&= - t^q \beta(t) \Bigl[ \mathcal{L}_{\epsilon(t)}(x(t)+\Gamma(t) \dot{x}(t), \lambda_t) - \mathcal{L}_{\epsilon(t)}(x_t, \lambda_t) \Bigr] - \frac{1}{2} t^q \epsilon(t) \|x(t) - x_t\|^2 \\
&\quad - t^q \epsilon(t) \Gamma(t) \langle x(t) - x_t, \dot{x}(t) \rangle - \frac{1}{2} t^q \epsilon(t) \Gamma(t)^2 \|\dot{x}(t)\|^2\\
&\quad +t^q\beta(t)\langle \Gamma(t)\dot{x}(t), \nabla_x \mathcal{L}_{\epsilon(t)}(x(t)+\Gamma(t) \dot{x}(t), \lambda_t) \rangle .
\end{split} \end{equation*}
Substituting the above inequality into \eqref{dE2}, we obtain
\begin{equation} \begin{split}
\dot{\mathcal{E}}_2(t) 
&\leq \bigl(1 + q t^{q-1} - \alpha\bigr) \langle x(t) - x_t, \dot{x}(t) \rangle + t^q \bigl(1 + q t^{q-1} - \alpha\bigr) \|\dot{x}(t)\|^2\ \\
&\quad- \langle x(t) - x_t, \dot{x}_t \rangle - t^q \langle \dot{x}(t), \dot{x}_t \rangle - t^{2q} \beta(t) \langle \dot{x}(t), \nabla_x \mathcal{L}_{\epsilon(t)}(x(t)+\Gamma(t) \dot{x}(t), \lambda_t) \rangle\\
&\quad +t^q\beta(t)\langle \Gamma(t)\dot{x}(t), \nabla_x \mathcal{L}_{\epsilon(t)}(x(t)+\Gamma(t) \dot{x}(t), \lambda_t) \rangle \\
&\quad- t^{q} \beta(t) \langle x(t) - x_t, A^\top(\lambda(t) - \lambda_t) \rangle - t^{2q} \beta(t) \langle \dot{x}(t), A^\top(\lambda(t) - \lambda_t) \rangle\\
&\quad - t^q \beta(t) \Bigl[ \mathcal{L}_{\epsilon(t)}(x(t)+\Gamma(t) \dot{x}(t), \lambda_t) - \mathcal{L}_{\epsilon(t)}(x_t, \lambda_t) \Bigr] - \frac{1}{2} t^q \epsilon(t) \|x(t) - x_t\|^2 \\
&\quad - t^q \epsilon(t) \Gamma(t) \langle x(t) - x_t, \dot{x}(t) \rangle - \frac{1}{2} t^q \epsilon(t) \Gamma(t)^2 \|\dot{x}(t)\|^2.
\end{split} \end{equation}

We now compute the derivatives of the remaining two components of the Lyapunov functional.
\begin{equation}
\begin{split}
\dot{\mathcal{E}}_3(t) &= \langle \lambda(t) - \lambda_t, \dot{\lambda}(t) - \dot{\lambda}_t \rangle ,
\end{split}
\end{equation}
\begin{equation}
\begin{split}
\dot{\mathcal{E}}_4(t) &= \frac{q(1-q)}{2} t^{q-2} \|x(t) - x_t\|^2 + \bigl(\alpha - 1 - q t^{q-1}\bigr) \langle x(t) - x_t, \dot{x}(t) - \dot{x}_t \rangle.
\end{split}
\end{equation}
Summing up the above estimates for
\(
\dot{\mathcal{E}}_1(t),\dot{\mathcal{E}}_2(t),\dot{\mathcal{E}}_3(t) \ \text{and} \ \dot{\mathcal{E}}_4(t),
\)
we obtain a preliminary differential inequality for the Lyapunov functional.

We next deal with all $\langle \nabla_x \mathcal{L}_{\epsilon(t)}(x(t)+\Gamma(t) \dot{x}(t), \lambda_t), \dot{x}(t) \rangle$ terms:
\begin{equation}
\begin{aligned}
&\left[ t^{2q} \beta(t) \left(1 + \dot{\Gamma}(t) - \frac{\alpha}{t^q} \Gamma(t)\right) - t^{2q} \beta(t) +t^q\beta(t)\Gamma(t)\right] 
\langle \nabla_x \mathcal{L}_{\epsilon(t)}(x(t)+\Gamma(t) \dot{x}(t), \lambda_t), \dot{x}(t) \rangle \\
&=\left[ t^{2q} \beta(t) \left(1 + \dot{\Gamma}(t) - \frac{\alpha}{t^q} \Gamma(t)\right) - t^{2q} \beta(t) +t^q\beta(t)\Gamma(t)\right]\\
& \quad \langle \nabla_x \mathcal{L}_{\epsilon(t)}(x(t)+\Gamma(t) \dot{x}(t), \lambda(t))-A^\top(\lambda-\lambda_t), \dot{x}(t) \rangle
\\
&=K(t)\langle \nabla_x \mathcal{L}_{\epsilon(t)}(x(t)+\Gamma(t) \dot{x}(t), \lambda(t)), \dot{x}(t) \rangle-K(t)\langle \lambda-\lambda_t, A\dot{x}(t) \rangle\\
& \leq |K(t)| \ \left[ \frac{\|\nabla_x \mathcal{L}_{\epsilon(t)}\bigl(x(t)+\Gamma(t) \dot{x}(t), \lambda(t)\bigr)\|^2}{2\eta_1(t)} + \frac{\eta_1(t)}{2} \|\dot{x}(t)\|^2 \right]-K(t)\langle \lambda-\lambda_t, A\dot{x}(t) \rangle,
\end{aligned}
\end{equation}
 where \(K(t) = t^{2q} \beta(t) \left(1 + \dot{\Gamma}(t) - \frac{\alpha-1}{t^q} \Gamma(t)\right) - t^{2q} \beta(t)=-t^{2q}\beta(t)\big(\frac{rq}{t^{q+1}}+\frac{\alpha-1}{t^q}(r_0+\tfrac{r}{t^q})\big)\).

We next focus on the terms containing $ \langle\cdot,\lambda(t) - \lambda_t \rangle$, an easy computation shows that
\begin{equation*}
\begin{split}
&t^{2q} \beta(t)^2 \Gamma(t) \big\langle A \nabla_x \mathcal{L}_{\epsilon(t)}\big(x(t)+\Gamma(t) \dot{x}(t), \lambda(t)\big), \lambda(t) - \lambda_t \big\rangle - t^q \beta(t) \big\langle A(x(t) - x_t), \lambda(t) - \lambda_t \big\rangle\\
&- t^{2q} \beta(t) \big\langle A \dot{x}(t), \lambda(t) - \lambda_t \big\rangle -K(t)\big\langle A\dot{x}(t),\lambda-\lambda_t \big\rangle
+ \big\langle \lambda(t) - \lambda_t, \dot{\lambda}(t) \big\rangle \\
&= t^{2q} \beta(t)^2 \Gamma(t) \big\langle A \nabla_x \mathcal{L}_{\epsilon(t)}\big(x(t)+\Gamma(t) \dot{x}(t), \lambda(t)\big), \lambda(t) - \lambda_t \big\rangle \\
&\quad - t^q \beta(t) \big\langle A\big(x(t) - x_t + t^q \dot{x}(t)\big), \lambda(t) - \lambda_t \big\rangle
-K(t)\big\langle A\dot{x}(t),\lambda-\lambda_t \big\rangle \\
&\quad - \big\langle  t^{2q} \beta(t)^2 \Gamma(t) A \nabla_x \mathcal{L}_{\epsilon(t)}\big(x(t)+\Gamma(t) \dot{x}(t), \lambda(t)\big),\lambda(t) - \lambda_t\big\rangle\\
&\quad + K(t)\big\langle A\dot{x}(t),\lambda-\lambda_t \big\rangle+t^q \beta(t) \big\langle \nabla_\lambda \mathcal{L}_{\epsilon(t)}\big(x(t) +t^q \dot{x}(t),\lambda(t)\big),\lambda(t) - \lambda_t\big\rangle \\
&=- t^q \beta(t) \big\langle A\big(x(t) - x_t + t^q \dot{x}(t)\big), \lambda(t) - \lambda_t \big\rangle
+ t^q \beta(t) \big\langle \nabla_\lambda \mathcal{L}_{\epsilon(t)}\big(x(t) +t^q \dot{x}(t),\lambda(t)\big),\lambda(t) - \lambda_t\big\rangle  \\
&=- t^q \beta(t) \big\langle A\big(x(t) - x_t + t^q \dot{x}(t)\big), \lambda(t) - \lambda_t \big\rangle \\
&\quad + t^q \beta(t) \big\langle A\big(x(t) + t^q \dot{x}(t)\big)-b-\frac{\epsilon(t)}{\beta(t)}\lambda(t),\lambda(t) - \lambda_t\big\rangle \\
&=- t^q \beta(t) \big\langle A\big(x(t) - x_t + t^q \dot{x}(t)\big), \lambda(t) - \lambda_t \big\rangle \\
&\quad + t^q \beta(t) \big\langle A\big(x(t) -x_t+ t^q \dot{x}(t)\big),\lambda(t) - \lambda_t\big\rangle
+ t^q \beta(t) \big\langle Ax_t-b-\frac{\epsilon(t)}{\beta(t)}\lambda(t),\lambda(t) - \lambda_t\big\rangle \\
&= t^q \beta(t) \big\langle \nabla_{\lambda} \mathcal{L}_{\epsilon(t)}\big(x_t,\lambda(t)\big),\lambda(t) - \lambda_t\big\rangle\\
&\leq  -\frac{1}{2} t^{q} \epsilon(t)\|\lambda(t) - \lambda_t\|^2,
\end{split}
\end{equation*}
where the last inequality uses
the \(\frac{\epsilon(t)}{\beta(t)}\)-strong convexity of \(-\mathcal{L}_{\epsilon(t)}(x_t,\cdot)\),  

\begin{equation*}
\begin{split}
\langle  \nabla_\lambda \mathcal{L}_{\epsilon(t)}(x_t, \lambda(t)) ,\lambda(t) - \lambda_t\rangle &\leq\mathcal{L}_{\epsilon(t)}(x_t, \lambda(t))-\mathcal{L}_{\epsilon(t)}(x_t, \lambda_t)-\frac{\epsilon(t)}{2\beta(t)} \|\lambda(t) - \lambda_t\|^2 \\
&\leq  -\frac{\epsilon(t)}{2\beta(t)} \|\lambda(t) - \lambda_t\|^2.
\end{split}
\end{equation*}

As $\Gamma(t)\to r_0$ as $t\to+\infty$, $\Gamma(t)$ becomes positive for sufficiently large $t$ provided $r_0>0$. When $r_0=0$ and $r\ge0$,  $\Gamma(t)\ge0$ for all $t\ge t_0$. Applying Young's inequality, the estimates below hold for sufficiently large $t$:
\begin{equation} \begin{split}
- t^q \epsilon(t) \Gamma(t) \langle x(t) - x_t, \dot{x}(t) \rangle
\leq t^q \epsilon(t) \Gamma(t) \ \frac{\|x(t) - x_t\|^2}{2\eta_0(t)}
+ t^q \epsilon(t) \Gamma(t) \ \frac{\eta_0(t)}{2} \|\dot{x}(t)\|^2,
\end{split} \end{equation}

\begin{equation} \begin{split}
t^{2q} \beta(t) \langle A(x(t) - x_t), \dot{\lambda}_t \rangle
&\leq \frac{t^{2q} \beta(t) \|A\|^2}{2 \eta_2(t)} \|x(t) - x_t\|^2 + \frac{t^{2q} \beta(t) \eta_2(t)}{2} \|\dot{\lambda}_t\|^2,
\end{split} \end{equation}
\begin{equation} \begin{split}
& \left( -\alpha + q t^{q-1} \right) \langle x(t) - x_t, \dot{x}_t \rangle \leq \left| -\alpha + q t^{q-1} \right| \left[ \frac{\|x(t) - x_t\|^2}{2 \eta_3(t)} + \frac{\eta_3(t)}{2} \|\dot{x}_t\|^2 \right] \\
&= \frac{\alpha - q t^{q-1}}{2 \eta_3(t)} \|x(t) - x_t\|^2 + \frac{\eta_3(t) (\alpha - q t^{q-1})}{2} \|\dot{x}_t\|^2,
\end{split} \end{equation}
\begin{equation} \begin{split}
- t^q \langle \dot{x}(t), \dot{x}_t \rangle \leq \frac{t^q}{2\eta_4(t)} \|\dot{x}(t)\|^2 + \frac{\eta_4(t) t^q}{2} \|\dot{x}_t\|^2,
\end{split} \end{equation}
\begin{equation} \begin{split} 
t^{2q} \beta(t) \Gamma(t) \langle A \dot{x}(t), \dot{\lambda}_t \rangle
&\leq \frac{t^{2q} \beta(t) \Gamma(t) \|A\|^2}{2\eta_5(t)} \|\dot{x}(t)\|^2 + \frac{\eta_5(t) t^{2q} \beta(t) \Gamma(t)}{2} \|\dot{\lambda}_t\|^2,
\end{split} \end{equation} 
\begin{equation} \begin{split} 
- \langle \lambda(t) - \lambda_t, \dot{\lambda}_t \rangle \leq \frac{1}{2\eta_6(t)} \|\lambda(t) - \lambda_t\|^2 + \frac{\eta_6(t)}{2} \|\dot{\lambda}_t\|^2.
\end{split} \end{equation}
Collecting the above estimates and applying Young's inequality to absorb mixed terms, we obtain
\begin{equation} \label{inequ:dEdt}\begin{split}
\dot{\mathcal{E}}(t) \leq{} & \bigl[t^{2q} \dot{\beta}(t) + 2q t^{2q-1} \beta(t) - t^{q} \beta(t)\bigr]\bigl(\mathcal{L}_{\epsilon(t)}(x(t)+\Gamma(t) \dot{x}(t), \lambda_t) - \mathcal{L}_{\epsilon(t)}(x_t, \lambda_t)\bigr) \\
&+ \left[ \frac{q(1-q)}{2} t^{q-2} - \frac{1}{2} t^q \epsilon(t) + \frac{t^q \epsilon(t) \Gamma(t)}{2\eta_0(t)} + \frac{\|A\|^2}{2\eta_2(t)} t^{2q} \beta(t) + \frac{\alpha - q t^{q-1}}{2\eta_3(t)} \right] \|x(t) - x_t\|^2 \\
&+ \Bigg[ t^q \bigl(1+q t^{q-1}-\alpha\bigr) - \frac{1}{2} t^q \epsilon(t) \Gamma(t)^2 + \frac{\eta_0(t)}{2} t^q \epsilon(t) \Gamma(t) + \frac{t^q}{2\eta_4(t)}\\
&\quad + \frac{\eta_1(t)}{2} 
\Big| t^{2q} \beta(t)\bigl(1+\dot{\Gamma}(t)-\frac{\alpha-1}{t^q}\Gamma(t)\bigr) - t^{2q} \beta(t)\Big|   + \frac{\|A\|^2}{2\eta_5(t)} t^{2q} \beta(t) \Gamma(t) \Bigg] \|\dot{x}(t)\|^2 \\
&+ \left[ \frac{\eta_3(t) (\alpha - q t^{q-1})}{2} + \frac{\eta_4(t) t^q}{2} \right] \|\dot{x}_t\|^2 + \left[ \frac{\eta_2(t)}{2} t^{2q} \beta(t) + \frac{\eta_5(t)}{2} \Gamma(t) t^{2q} \beta(t) + \frac{\eta_6(t)}{2} \right] \|\dot{\lambda}_t\|^2 \\
&+ \left( -t^{2q} \beta(t)^2 \Gamma(t) + \frac{1}{2\eta_1(t)} \Big| t^{2q} \beta(t)\bigl(1-\frac{rq}{t^{q+1}} - \frac{\alpha-1}{t^q}\Gamma(t))- t^{2q} \beta(t)\Big| \right) \\
&\quad \|\nabla_x \mathcal{L}_{\epsilon(t)}(x(t)+\Gamma(t) \dot{x}(t), \lambda(t))\|^2 \\
&+ \left[ -\frac{1}{2}t^q \epsilon(t) + \frac{1}{2\eta_6(t)}  - t^{2q} \beta(t)^2\Gamma(t) \|A\|^2 + t^{2q} \beta(t)^2\Gamma(t) \|A\|^2\right] \|\lambda(t) - \lambda_t\|^2 \\
&+ t^{2q} \beta(t)\left(\frac{\epsilon(t)}{2\beta(t)}\right)' (\|x(t)+\Gamma(t)\dot{x}(t)\|^2-\|x_t\|^2).\\
\end{split} \end{equation}

From the definition of the Lyapunov function\eqref{eq:energy} ,we obtain
\begin{equation}\label{inequ:E}\begin{split}
\mathcal{E}(t) &\leq t^{2q} \beta(t) \bigl(\mathcal{L}_{\epsilon(t)}\bigl(x(t) + \Gamma(t)\dot{x}(t), \lambda_t\bigr) - \mathcal{L}_{\epsilon(t)}(x_t, \lambda_t)\bigr) + \left( \frac{\alpha-1-q t^{q-1}}{2} + 1 \right) \|x(t)-x_t\|^2\\
&\quad + t^{2q} \|\dot{x}(t)\|^2 + \frac{1}{2} \|\lambda(t) - \lambda_t\|^2.
\end{split} \end{equation}

Combining the above inequalities \eqref{inequ:dEdt} and \eqref{inequ:E} ,we obtain the following estimate:
\begin{equation} \begin{split}
\dot{\mathcal{E}}(t) &+ \frac{W}{t^m} \mathcal{E}(t) \\
\leq{} & \bigl[t^{2q} \dot{\beta}(t) + \bigl(\frac{2q}{t} + \frac{W}{t^m} - \frac{1}{t^{q}}\bigr) t^{2q} \beta(t)\bigr]\bigl(\mathcal{L}_{\epsilon(t)}(x(t)+\Gamma(t) \dot{x}(t), \lambda_t)- \mathcal{L}_{\epsilon(t)}(x_t, \lambda_t)\bigr) \\
&+ \left[ \frac{q(1-q)}{2} t^{q-2} - \frac{1}{2} t^q \epsilon(t) + \frac{t^q \epsilon(t) \Gamma(t)}{2\eta_0(t)} + \frac{\|A\|^2}{2\eta_2(t)} t^{2q} \beta(t) + \frac{\alpha - q t^{q-1}}{2\eta_3(t)} + \frac{W}{t^m} \ \frac{\alpha+1-q t^{q-1}}{2} \right] \\
&\quad\|x(t)-x_t\|^2 \\
&+ \Bigg[ t^q \bigl(1+q t^{q-1}-\alpha\bigr) - \frac{1}{2} t^q \epsilon(t) \Gamma(t)^2 + \frac{\eta_0(t)}{2} t^q \epsilon(t) \Gamma(t) \\
&\quad + \frac{\eta_1(t)}{2} \ \Big| t^{2q} \beta(t)\bigl(\frac{rq}{t^{q+1}} + \frac{\alpha-1}{t^q} \Gamma(t) \bigl) \Big| + \frac{t^q}{2\eta_4(t)} + \frac{\|A\|^2}{2\eta_5(t)} t^{2q} \beta(t) \Gamma(t) + \frac{W}{t^m} \ t^{2q} \Bigg] \|\dot{x}(t)\|^2 \\
&+ \left[ \frac{\eta_3(t) (\alpha - q t^{q-1})}{2} + \frac{\eta_4(t) t^q}{2} \right] \|\dot{x}_t\|^2 + \left[ t^{2q} \beta(t) \bigl(\frac{\eta_2(t)}{2} + \frac{\eta_5(t)}{2} \Gamma(t)\bigr) + \frac{\eta_6(t)}{2} \right] \|\dot{\lambda}_t\|^2\\
&+ \left[ -t^{2q} \beta(t)^2 \Gamma(t)+ \frac{1}{2\eta_1(t)} \Big| t^{2q} \beta(t)\bigl(\frac{rq}{t^{q+1}} + \frac{\alpha-1}{t^q}\Gamma(t)\bigr) \Big| \right] \|\nabla_x \mathcal{L}_{\epsilon(t)}(x(t)+\Gamma(t) \dot{x}(t), \lambda(t))\|^2 \\
&+ \left[ -\frac{1}{2}t^q \epsilon(t) + \frac{1}{2\eta_6(t)} + \frac{W}{2 t^m} \right] \|\lambda(t)-\lambda_t\|^2\\
&+ t^{2q} \beta(t)\left(\frac{\epsilon(t)}{2\beta(t)}\right)' (\|x(t)+\Gamma(t)\dot{x}(t)\|^2-\|x_t\|^2).
\end{split} \end{equation}

We now analyze the coefficients of the various terms in the time derivative of the Lyapunov functional.

Consider the coefficient of $\mathcal{L}_{\epsilon(t)}(x(t)+\Gamma(t)\dot{x}(t),\lambda_t)-\mathcal{L}_{\epsilon(t)}(x_t,\lambda_t)$. There exist $t_1 \ge t_0$ and $P_1 \ge 0$ such that for all $t \ge t_1$,
\[
\begin{split}
t^{2q} \dot{\beta}(t) + \left(\frac{2q}{t} + \frac{W}{t^m} - \frac{1}{t^q}\right) t^{2q} \beta(t)
\le -P_1 t^{q} \beta(t) .
\end{split}
\]
Next, we examine the coefficient of $\|x(t) - x_t\|^2$. Choosing $\eta_0(t) =4(r_0+1)$, $\eta_2(t) = t$, and letting $\eta_3(t) = \frac{4\alpha}{t^q \epsilon(t)}$, we consider two cases.\\
If $r_0 = 0$ and $r \ge 0$, the coefficient reduces to
\[
\frac{q(1-q)}{2} t^{q-2}
- \frac{1}{2} t^q \epsilon(t)
+ \frac{r \epsilon(t)}{2\eta_0(t)}
+ \frac{\|A\|^2}{2\eta_2(t)} t^{2q} \beta(t)
+ \frac{\alpha - q t^{q-1}}{2\eta_3(t)}
+ \frac{W}{t^m} \frac{\alpha + 1 - q t^{q-1}}{2},
\]
If $r_0 > 0$ and $r \in \mathbb{R}$, the coefficient becomes
\[
\frac{q(1-q)}{2} t^{q-2}
+ \frac{1}{2} t^q \epsilon(t) \left( \frac{\Gamma(t)}{\eta_0(t)} - 1 \right)
+ \frac{\|A\|^2}{2\eta_2(t)} t^{2q} \beta(t)
+ \frac{\alpha - q t^{q-1}}{2\eta_3(t)}
+ \frac{W}{t^m} \frac{\alpha + 1 - q t^{q-1}}{2},
\]
Under the assumption \eqref{H0}, there exist $t_2 \ge t_0$ and $P_2 > 0$ such that for all $t \ge t_2$,
\[
\begin{aligned}
& \left[
    \frac{q(1-q)}{2} t^{q-2}
    - \frac{1}{2} t^q \epsilon(t)
    + \frac{\epsilon(t) t^q}{2\eta_0(t)} \left( r_0 + \frac{r}{t^q} \right)
    + \frac{\|A\|^2}{2\eta_2(t)} t^{2q} \beta(t)
    + \frac{\alpha - q t^{q-1}}{2\eta_3(t)}
    + \frac{W(\alpha+1-q t^{q-1})}{2 t^m}
\right] \\
& \|x(t) - x_t\|^2 \le - P_2 t^q \epsilon(t) \|x(t) - x_t\|^2.
\end{aligned}
\]

We now turn to the coefficient of $\|\dot{x}(t)\|^2$. Choosing $\eta_1(t) = \frac{1}{4(r_0+1)\beta(t)}$ , $\eta_4(t) = \frac{1}{\alpha-1}$, $\eta_5(t) = t$, and substituting $\eta_0(t)=4(r_0+1)$, we again consider two cases.\\
If $r_0 = 0$ and $r \ge 0$, the coefficient becomes
\begin{equation*}
\begin{split}
&t^q \left( 1 + q t^{q-1} - \alpha \right)
    -  \frac{r^2 \epsilon(t)}{2t^q}
    + \frac{r\eta_0(t)}{2}  \epsilon(t)
    + \frac{\eta_1(t)}{2} t^{2q} \beta(t) \left( \frac{r q}{t^{q+1}} + \frac{(\alpha-1) r}{t^{2q}} \right) 
    + \frac{t^q}{2\eta_4(t)} \\
    & + \frac{\|A\|^2}{2\eta_5(t)} t^q \beta(t) r
    + \frac{W}{t^m} t^{2q} .
\end{split}
\end{equation*}
As $t\to\infty$, this expression is asymptotically equivalent to
    \[
    \begin{split}
     t^q (1-\alpha)
    - \frac{r^2 \epsilon(t)}{t^q}
    + \frac{r \eta_0(t)}{2} \epsilon(t)
    + \frac{(\alpha-1) r}{2} \eta_1(t) \beta(t)
    + \frac{1}{2\eta_4(t)} t^q
    + \frac{\|A\|^2}{2\eta_5(t)} r t^q \beta(t)
    + \frac{W}{t^m} t^{2q},
    \end{split}
    \]
If $r_0 > 0$ and $r \in \mathbb{R}$, the coefficient becomes
\begin{equation*}
\begin{split}
     &t^q \left( 1 + q t^{q-1} - \alpha \right)
    - \frac{r_0}{2} t^q \epsilon(t)
   \left(r_0-\eta_0(t)\right)
    + \frac{\eta_1(t)}{2} t^{2q} \beta(t) \left( \frac{r q}{t^{q+1}} + \frac{(\alpha-1) r_0}{t^q} \right) 
    + \frac{t^q}{2\eta_4(t)} \\
   & + \frac{\|A\|^2}{2\eta_5(t)} r_0 t^{2q} \beta(t)
    + \frac{W}{t^m} t^{2q},
\end{split}
\end{equation*}
    As $t\to\infty$, this expression behaves asymptotically as
    \[
    \begin{split}
    t^q (1-\alpha)
    - \frac{r_0^2}{2} t^q \epsilon(t)
    + \frac{r_0}{2} \eta_0(t) t^q \epsilon(t)
    + \frac{(\alpha-1) r_0}{2} \eta_1(t) t^q \beta(t)
    + \frac{1}{2\eta_4(t)} t^q
    + \frac{\|A\|^2}{2\eta_5(t)} r_0 t^{2q} \beta(t)
    + \frac{W}{t^m} t^{2q}.
    \end{split}
    \]
Under the assumption $m > q$, $t^{2q-m} < t^q$. Also under the assumption $\lim_{t\to+\infty}{\epsilon(t)} = 0,$ 
and $\lim_{t\to+\infty} \frac{t^{q-1}\beta(t)}{\epsilon(t)} = 0,$ hence, there exist $t_3 \ge t_0$ and $\lim_{t\to+\infty} \frac{t^{q-1}\beta(t)}{\epsilon(t)} = 0,$ $P_3 > 0$ such that for all $t \ge t_3$,
\[
\begin{aligned}
&\Bigg[
t^q\big(1 + q t^{q-1} - \alpha\big)
- \frac{1}{2} t^q \epsilon(t) \left(r_0 + \frac{r}{t^q}\right)^2 + \frac{\eta_0(t)}{2} t^q \epsilon(t) \left(r_0 + \frac{r}{t^q}\right)+ \frac{t^q}{2\eta_4(t)}\\ 
&\quad+ \frac{\eta_1(t)}{2} \left| t^{2q} \beta(t) \left( \frac{r q}{t^{q+1}} + \frac{\alpha-1}{t^q}\left(r_0 + \frac{r}{t^q}\right) \right) \right| 
+ \frac{\|A\|^2}{2\eta_5(t)} t^{2q} \beta(t)\left(r_0 + \frac{r}{t^q}\right)
+ \frac{W}{t^m} t^{2q}
\Bigg] \|\dot{x}(t)\|^2 \\
&\le - P_3 t^q \|\dot{x}(t)\|^2.
\end{aligned}
\]

Next, we consider the coefficient of $\|\dot{x}_t\|^2$. Substituting $\eta_3(t) = \frac{4\alpha}{t^q \epsilon(t)}$ and $\eta_4(t) = \frac{1}{\alpha-1}$, the coefficient behaves asymptotically as
\(
\frac{1}{2(\alpha-1)} t^q + \frac{\alpha}{2} \eta_3(t).
\)
Thus, for all $t \ge \max\{t_2, t_3\}$, 
\begin{equation*}
\begin{split}
\left( \frac{\eta_3(t) (\alpha - q t^{q-1})}{2} + \frac{\eta_4(t)}{2} t^q \right) \|\dot{x}_t\|^2 \le \left( \frac{\alpha}{2} \eta_3(t) + \frac{1}{2(\alpha-1)} t^q \right) \|\dot{x}_t\|^2  \le \left( \frac{2\alpha^2}{ t^q \epsilon(t)} + \frac{t^q}{\alpha-1} \right) \|\dot{x}_t\|^2.
\end{split}
\end{equation*}

We now analyze the coefficient of $\|\dot{\lambda}_t\|^2$. Choosing $\eta_2(t) = t$, $\eta_5(t) = t$, and $\eta_6(t) = \frac{2}{t^q \epsilon(t)}$, we consider two cases.\\
If $r_0 = 0$ and $r \ge 0$, the coefficient becomes
    \[
    \frac{1}{2} t^{2q+1} \beta(t) + \frac{r}{2} t^{q+1} \beta(t) + \frac{1}{t^q \epsilon(t)},
    \]
If $r_0 > 0$ and $r \in \mathbb{R}$, the coefficient becomes
    \[
    \frac{1}{2} t^{2q+1} \beta(t) + \frac{r_0}{2} t^{2q+1} \beta(t) + \frac{r}{2} t^{q+1} \beta(t) + \frac{1}{t^q \epsilon(t)}.
    \]
Hence, there exists $t \ge \max\{t_1, t_2, t_3\}$ such that
\begin{equation*}
\begin{aligned}
\Bigg[ t^{2q} \beta(t) \bigg( \frac{\eta_2(t)}{2} + \frac{\eta_5(t)}{2} \Gamma(t) \bigg) + \frac{\eta_6(t)}{2} \Bigg] \|\dot{\lambda}_t\|^2 
\le \bigg( \frac{r_0+1}{2} t^{2q+1} \beta(t) + \frac{r}{2} t^{q+1} \beta(t) + \frac{1}{t^q \epsilon(t)} \bigg) \|\dot{\lambda}_t\|^2.
\end{aligned}
\end{equation*}

We now examine the coefficient of $\|\lambda(t) - \lambda_t\|^2$. Substituting $\eta_6(t) = \frac{2}{t^q \epsilon(t)}$ and assuming $\textup{(H}_0\textup{)}$ , and using $\epsilon(t) > \frac{\alpha}{t^{q+m}}$ for all sufficiently large $t$, so there exist $t_4 \ge t_0$ and $P_4 > 0$ such that for all $t \ge t_4$,
\[
\left( \frac{1}{2\eta_6(t)} - \frac{1}{2}t^q \epsilon(t) + \frac{W}{2t^m} \right) \|\lambda(t) - \lambda_t\|^2
\le - P_4 t^q \epsilon(t) \|\lambda(t) - \lambda_t\|^2.
\]

Finally, we consider the coefficient of $\|\nabla_x \mathcal{L}_{\epsilon(t)}(x(t) + \Gamma(t) \dot{x}(t), \lambda(t))\|^2$. Substituting $\eta_1(t) = \frac{1}{4(r_0+1)\beta(t)} $, we consider two cases.\\
If $r_0 = 0$ and $r > 0$, the coefficient becomes
    \[
    - t^{2q} \beta(t)^2 \ \frac{r}{t^q} + \frac{1}{2\eta_1(t)} \left[ t^{2q} \beta(t) \left( \frac{r q}{t^{q+1}} + \frac{(\alpha-1) r}{t^{2q}} \right) \right],
    \]
As $t\to\infty$, this expression behaves asymptotically as
    \(
    - r t^q \beta(t)^2.
    \)
When $r=0$, the terms involving $\Gamma(t)$ vanish, and the corresponding estimate remains valid with $P_
5=0$.\\
If $r_0 > 0$ and $r \in \mathbb{R}$, the coefficient becomes
    \[
    - t^{2q} \beta(t)^2 r_0 + \frac{1}{2\eta_1(t)} t^{2q} \beta(t) \cdot \frac{(\alpha-1) r_0}{t^q}
    = - r_0 t^{2q} \beta(t)^2 + {2r_0 (r_0+1) (\alpha-1)} t^q \beta(t)^2,
    \]
As $t\to\infty$, this expression is asymptotically equivalent to 
    \(
    - r_0 t^{2q} \beta(t)^2.
    \)
Thus, there exist $t_5 \ge t_0$ and $P_5 \ge 0$ such that for all $t \ge t_5$,
\begin{equation*}
\begin{split}
&\left[
- t^{2q} \beta(t)^2 \left( r_0 + \frac{r}{t^q} \right)
+ \frac{1}{2\eta_1(t)} t^{2q} \beta(t) \left( \frac{r q}{t^{q+1}} + \frac{\alpha-1}{t^q} \left( r_0 + \frac{r}{t^q} \right) \right)
\right]\\
&\|\nabla_x \mathcal{L}_{\epsilon(t)}(x(t)+\Gamma(t) \dot{x}(t), \lambda(t))\|^2 \\
&\quad \le - P_5 t^{q} \beta(t)^2 \|\nabla_x \mathcal{L}_{\epsilon(t)}(x(t)+\Gamma(t) \dot{x}(t), \lambda(t))\|^2.
\end{split}
\end{equation*}

Collecting the above estimates, we finally obtain,
\begin{equation} \label{inequ:E+dE} \begin{split}
\dot{\mathcal{E}}(t) + \frac{W}{t^m} \mathcal{E}(t)
\le{} & - P_1 t^{q} \beta(t) \bigl(\mathcal{L}_{\epsilon(t)}(x(t)+\Gamma(t) \dot{x}(t), \lambda_t) - \mathcal{L}_{\epsilon(t)}(x_t, \lambda_t)\bigr)
- P_2 \epsilon(t) t^q \|x(t)-x_t\|^2\\
&- P_3 t^q \|\dot{x}(t)\|^2 
- P_4 \epsilon(t) t^q \|\lambda(t) - \lambda_t\|^2 \\
&- P_5 t^{q} \beta(t)^2 \|\nabla_x \mathcal{L}_{\epsilon(t)}(x(t)+\Gamma(t) \dot{x}(t), \lambda(t))\|^2\\
&+ \left( \frac{2\alpha^2}{ t^q \epsilon(t)} + \frac{t^q}{\alpha-1} \right) \|\dot{x}_t\|^2+ \left( \frac{r_0+1}{2} t^{2q+1} \beta(t) + \frac{r}{2} t^{q+1} \beta(t) + \frac{1}{t^q \epsilon(t)} \right) \|\dot{\lambda}_t\|^2\\
&+ t^{2q} \beta(t)\left(\frac{\epsilon(t)}{2\beta(t)}\right)' (\|x(t)+\Gamma(t) \dot{x}(t)\|^2-\|x_t\|^2). 
\end{split} \end{equation}

\end{proof}

\textbf{Remark 3.1} The assumptions in \eqref{H0} are nonempty and are satisfied by standard power-type coefficients. For example, let
$\beta(t)=t^{s},\ \epsilon(t)=t^{-p},$
where
$0<q<m<1,\ s>0,\ p>0,$
and suppose that
$q+s+p<1,\ p<q+m.$
Then
$t^{q-1}\frac{\beta(t)}{\epsilon(t)}=t^{q-1+s+p}\longrightarrow 0,$
and
$\frac{\epsilon(t)}{\alpha t^{-(q+m)}}=\frac{1}{\alpha}t^{q+m-p}\longrightarrow +\infty,$
so that $\epsilon(t)>\frac{\alpha}{t^{q+m}}$ for all sufficiently large $t$. Moreover,
$t^{q}\frac{\dot{\beta}(t)}{\beta(t)}=s t^{q-1}\longrightarrow 0<1.$
Hence all the conditions in \eqref{H0} are fulfilled.

The following theorem provides convergence estimates under two sufficient eventual comparison regimes for the perturbation terms appearing in \eqref{inequ:deal}.

\begin{theorem}\label{Th3.2}
Suppose that $0 < W < \min\left\{ \frac{\alpha}{2},\ \frac{\alpha}{2(\alpha+1)} \right\}$ and all assumptions of Theorem \ref{Th3.1} hold. And further assume that $\frac{\epsilon(t)}{\beta(t)} \in C^{2}\bigl([t_{0},+\infty);(0,+\infty)\bigr)$. Then the following assertions hold:

\begin{enumerate}
\item[(i)]
Suppose that $\exists\ T_1 \ge t_0,$ \text{and a constant }$\rho_1\in(0,W),$
such that for every $t \ge T_1$,\\
\[\left[\left(\frac{\epsilon(t)}{\beta(t)}\right)'^2\left(\frac{\beta(t)}{\epsilon(t)}\right)^2\right]' \ge -\frac{\rho_1}{t^m} \left(\frac{\epsilon(t)}{\beta(t)}\right)'^2\left(\frac{\beta(t)}{\epsilon(t)}\right)^2,
\]
\[
t^{2q+1} \beta(t)\left(\frac{\epsilon(t)}{\beta(t)}\right)'^2\left(\frac{\beta(t)}{\epsilon(t)}\right)^2 \ge  t^{2q} \beta(t)\left|\left(\frac{\epsilon(t)}{2\beta(t)}\right)' \right|,
\]
and
\[
\lim_{t\to\infty} t^{2q+1+m} \beta(t)\left(\frac{\epsilon(t)}{\beta(t)}\right)'^2\left(\frac{\beta(t)}{\epsilon(t)}\right)^2 = 0,
\]
then the following estimates hold:
\begin{equation*}
\begin{aligned}
\|x(t)-x_t\|^2 \leq O\!\left( t^{2q+m+1} \beta(t) \left(\frac{\epsilon(t)}{\beta(t)}\right)'^2\left(\frac{\beta(t)}{\epsilon(t)}\right)^2 \right), 
\end{aligned}
\end{equation*}
\begin{equation*}
\begin{aligned}
\|\lambda(t)-\lambda_t\|^2 \leq O\!\left( t^{2q+m+1} \beta(t) \left(\frac{\epsilon(t)}{\beta(t)}\right)'^2\left(\frac{\beta(t)}{\epsilon(t)}\right)^2 \right),
\end{aligned}
\end{equation*}
\begin{equation*}
\begin{aligned}
\|\dot{x}(t)\|^2 \le O\left( t^{m+1} \beta(t) \left(\frac{\epsilon(t)}{\beta(t)}\right)'^2\left(\frac{\beta(t)}{\epsilon(t)}\right)^2\right),
\end{aligned}
\end{equation*}
\begin{equation*}
\begin{aligned}
\|Ax(t)-b\| \leq O\left( \frac{\epsilon(t)}{\beta(t)} + t^{\frac{2q+m+1}{2}} \sqrt{\beta(t)}\left|\left(\frac{\epsilon(t)}{\beta(t)}\right)'\right| \frac{\beta(t)}{\epsilon(t)} \right),
\end{aligned}
\end{equation*}
\[
\mathcal{L}_{\epsilon(t)}(x(t)+\Gamma(t)\dot{x}(t), \lambda_t) - \mathcal{L}_{\epsilon(t)}(x_t, \lambda_t) \le O\left({t^{m+1}}\left(\frac{\epsilon(t)}{\beta(t)}\right)'^2\left(\frac{\beta(t)}{\epsilon(t)}\right)^2 \right),
\]

\begin{equation*}
\begin{split}
&|f(x(t)+\Gamma(t) \dot{x}(t)) - f(x^*)|\\
&\leq O\left( \frac{\epsilon(t)}{\beta(t)} + {t^{\frac{2q+m+1}{2}}} \left|\left(\frac{\epsilon(t)}{\beta(t)}\right)'\right| \frac{\beta(t)}{\epsilon(t)}\sqrt{\beta(t)}  +{t^{m+1}}\left(\frac{\epsilon(t)}{\beta(t)}\right)'^2\left(\frac{\beta(t)}{\epsilon(t)}\right)^2 \right).
\end{split}
\end{equation*}

\item[(ii)]
Suppose that $\exists\ T_2 \ge t_0, \text{and a constant }\rho_2\in(0,W),$ such that for every $t \ge T_2$,\\
\[
\left(-\frac{\epsilon(t)}{\beta(t)}\right)''
\ge -\frac{\rho_2}{t^m}
\left(-\frac{\epsilon(t)}{\beta(t)}\right)',
\quad
t^{2q}\beta(t)\left|\left(\frac{\epsilon(t)}{2\beta(t)}\right)'\right|\ge t^{2q+1} \beta(t)\left(\frac{\epsilon(t)}{\beta(t)}\right)'^2\left(\frac{\beta(t)}{\epsilon(t)}\right)^2,
\]
and
\[
\lim_{t\to\infty} t^{2q+m}\beta(t)\left|\left(\frac{\epsilon(t)}{2\beta(t)}\right)'\right| = 0,
\]
then the following estimates hold:
\begin{equation*}
\begin{aligned}
\|x(t)-x_t\|^2 &\leq O\left( t^{2q+m} \beta(t) \left(-\frac{\epsilon(t)}{\beta(t)}\right)' \right),\
\|\lambda(t)-\lambda_t\|^2 \leq O\left( t^{2q+m} \beta(t) \left(-\frac{\epsilon(t)}{\beta(t)}\right)' \right),
\end{aligned}
\end{equation*}
\begin{equation*}
\begin{aligned}
\|\dot{x}(t)\|^2 &\leq O\left( t^{m} \beta(t) \left(-\frac{\epsilon(t)}{\beta(t)}\right)' \right),\
\|Ax(t)-b\| \leq O\left( \frac{\epsilon(t)}{\beta(t)} + {t^{\frac{2q+m}{2}}} \sqrt{ \beta(t)\left(-\frac{\epsilon(t)}{\beta(t)}\right)'} \right),
\end{aligned}
\end{equation*}
\[
\mathcal{L}_{\epsilon(t)}(x(t)+\Gamma(t)\dot{x}(t), \lambda_t) - \mathcal{L}_{\epsilon(t)}(x_t, \lambda_t)\leq  O\left({t^{m}} \left(-\frac{\epsilon(t)}{\beta(t)}\right)'  \right),
\]
\[
|f(x(t)+\Gamma(t) \dot{x}(t)) - f(x^*)| \leq O\left( \frac{\epsilon(t)}{\beta(t)} + {t^{\frac{2q+m}{2}}} \sqrt{\beta(t)\left(-\frac{\epsilon(t)}{\beta(t)}\right)'} + {t^{m} } \left(-\frac{\epsilon(t)}{\beta(t)}\right)' \right).
\]
\end{enumerate}
In both cases, the entire primal-dual trajectory converges
strongly to the minimum-norm saddle point:
\[
\lim_{t\to+\infty}
\|(x(t),\lambda(t))-(x^*,\lambda^*)\|=0,
\]
\end{theorem}

\begin{proof}

We now derive the asymptotic convergence rates associated with the Lyapunov estimate established in Theorem \ref{Th3.1}. The proof relies on transforming the differential inequality satisfied by the Lyapunov functional into explicit decay estimates through a suitable integrating-factor argument. Let $y^*=(x^*,\lambda^*)$ and $y_t=(x_t,\lambda_t)$, then using Lemma \ref{lem2.1} we obtain
\begin{equation} \label{xy1} \begin{split}
\|x_t\|^2 \leq \|y_t\|^2 \leq \|y^*\|^2,
\end{split} \end{equation}
and for a.e. $t \ge t_0$
\begin{equation}  \label{xy2} \begin{split}
\max\left\{\|\dot{x}_t\|^2, \|\dot{\lambda}_t\|^2\right\} \leq \|\dot{y}_t\|^2 \leq \left(\frac{\epsilon(t)}{\beta(t)}\right)'^2\left(\frac{\beta(t)}{\epsilon(t)}\right)^2 \|y_t\|^2 \leq \left(\frac{\epsilon(t)}{\beta(t)}\right)'^2\left(\frac{\beta(t)}{\epsilon(t)}\right)^2 \|y^*\|^2.
\end{split} \end{equation}
By neglecting the nonpositive dissipation terms in \eqref{inequ:E+dE} and using \eqref{xy1} and \eqref{xy2}, we obtain the following simplified upper estimate valid for all $t \ge \max\{t_1, \dots, t_5\}$,
\begin{equation} \label{inequ:deal} \begin{split}
\dot{\mathcal{E}}(t) + \frac{W}{t^m} \mathcal{E}(t)
&\le \left( \frac{2\alpha^2}{ t^q \epsilon(t)} + \frac{t^q}{\alpha-1} \right) \left(\frac{\epsilon(t)}{\beta(t)}\right)'^2\left(\frac{\beta(t)}{\epsilon(t)}\right)^2\|y^*\|^2 - t^{2q} \beta(t) \left(\frac{\epsilon(t)}{2\beta(t)}\right)' \|x_t\|^2\\
&\quad + \left( \frac{r_0+1}{2} t^{2q+1} \beta(t) + \frac{r}{2} t^{q+1} \beta(t) + \frac{1}{t^{q} \epsilon(t)} \right) \left(\frac{\epsilon(t)}{\beta(t)}\right)'^2\left(\frac{\beta(t)}{\epsilon(t)}\right)^2\|y^*\|^2 .
\end{split} \end{equation}
Under assumption \eqref{H0}, $\epsilon(t)>\frac{\alpha}{t^{q+m}}$ for all sufficiently large $t$,  $t^{3q+1}\beta(t)\epsilon(t) \ge \alpha \beta(t) t^{2q+1-m}$. Since $0<q<m<1$ and $\beta \to +\infty$, it follows that $\lim_{t\to+\infty}t^{3q+1}\beta(t)\epsilon(t)=+\infty.$ Hence, by enlarging $T_0$ if necessary,
\(
\frac{1}{t^q\epsilon(t)}
\le
t^{2q+1}\beta(t).\) 
Moreover, the remaining terms of orders $t^q<t^{q+1} \beta(t) <t^{2q+1} \beta(t)$ for all sufficiently large $t$. Therefore, \eqref{inequ:deal} reduces to
\begin{equation} \label{two cases}
\begin{split}
\dot{\mathcal{E}}(t) + \frac{W}{t^m} \mathcal{E}(t)
\le C\Bigg[
& t^{2q+1}\beta(t)
\left(\frac{\epsilon(t)}{\beta(t)}\right)'^{2}
\left(\frac{\beta(t)}{\epsilon(t)}\right)^{2} + t^{2q}\beta(t)
\left|
\left(\frac{\epsilon(t)}{2\beta(t)}\right)'
\right|
\Bigg]\|y^*\|^2.
\end{split}
\end{equation}
We therefore consider the two sufficient eventual comparison regimes between the remaining perturbation terms appearing on the right-hand side of \eqref{two cases}.\\
\paragraph{Case (I):} \[t^{2q+1} \beta(t) \left(\frac{\epsilon(t)}{\beta(t)}\right)'^2\left(\frac{\beta(t)}{\epsilon(t)}\right)^2\ge \max\left\{ \frac{1}{t^{q} \epsilon(t)}  \left(\frac{\epsilon(t)}{\beta(t)}\right)'^2\left(\frac{\beta(t)}{\epsilon(t)}\right)^2,  t^{2q} \beta(t)\left|\left(\frac{\epsilon(t)}{2\beta(t)}\right)' \right|\right\}\] 
Consequently, the differential inequality reduces to the following asymptotic form
\begin{equation}
\dot{\mathcal{E}}(t) + \frac{W}{t^m} \mathcal{E}(t)
\leq C_0 t^{2q+1} \beta(t) \left(\frac{\epsilon(t)}{\beta(t)}\right)'^2\left(\frac{\beta(t)}{\epsilon(t)}\right)^2 \|y^*\|^2,
\end{equation}
that is
\[
\dot{\mathcal{E}}(t) + \frac{W}{t^m} \mathcal{E}(t) \le C_1 t^{2q+1}\left(\frac{\epsilon(t)}{\beta(t)}\right)'^2\left(\frac{\beta(t)}{\epsilon(t)}\right)^2 \beta(t).
\]

Since the differential inequality contains the damping term
\(
\frac{W}{t^m} \mathcal{E}(t),
\)
with $0 < m < 1$, it is natural to introduce an integrating factor in order to absorb the dissipative contribution. Define
\[
\phi(t) = \int_{t_0}^t \frac{W}{s^m} ds = \frac{W}{1-m} \left( t^{1-m} - t_0^{1-m} \right),
\]
where \ \(\phi'(t) = \frac{W}{t^m}\) exactly matches the dissipation structure.
Multiplying the inequality by \(e^{\phi(t)}\)  transforms the left-hand side into an exact derivative, this yields
\[
e^{\phi(t)} \left( \dot{\mathcal{E}}(t) + \frac{W}{t^m} \mathcal{E}(t) \right) \le C_1 e^{\phi(t)} t^{2q+1} \beta(t) \left(\frac{\epsilon(t)}{\beta(t)}\right)'^2\left(\frac{\beta(t)}{\epsilon(t)}\right)^2,
\]
that is
\begin{equation} \label{under_use_1}\begin{split}
\frac{d}{dt} \left( e^{\phi(t)} \mathcal{E}(t) \right) \le C_1 e^{\phi(t)} t^{2q+1} \beta(t)  \left(\frac{\epsilon(t)}{\beta(t)}\right)'^2\left(\frac{\beta(t)}{\epsilon(t)}\right)^2.
\end{split} \end{equation} 
To estimate the right-hand side, we introduce a comparison quantity whose derivative reproduces the same asymptotic structure as the dominant perturbation term. Let
$t_6:=\max\{T_0,T_1\},\
C_2:=\frac{W-\rho_1}{W}\in(0,1).$ A direct computation shows that for \(t_6 \ge t_0\)
\begin{equation} \label{under_use_2} \begin{split}
\frac{d}{dt} \left( e^{\phi(t)} t^{2q+1+m} \beta(t) \left(\frac{\epsilon(t)}{\beta(t)}\right)'^2\left(\frac{\beta(t)}{\epsilon(t)}\right)^2 \right)
= e^{\phi(t)} \ \frac{W}{t^m} t^{2q+1+m} \beta(t)  \left(\frac{\epsilon(t)}{\beta(t)}\right)'^2\left(\frac{\beta(t)}{\epsilon(t)}\right)^2\\
+ e^{\phi(t)} (2q+1+m) t^{2q+m} \beta(t)  \left(\frac{\epsilon(t)}{\beta(t)}\right)'^2\left(\frac{\beta(t)}{\epsilon(t)}\right)^2
+ e^{\phi(t)} t^{2q+1+m} \dot{\beta}(t)  \left(\frac{\epsilon(t)}{\beta(t)}\right)'^2\left(\frac{\beta(t)}{\epsilon(t)}\right)^2\\ +e^{\phi(t)}t^{2q+m+1} \beta(t)  \left(\left(\frac{\epsilon(t)}{\beta(t)}\right)'^2\left(\frac{\beta(t)}{\epsilon(t)}\right)^2 \right)'
\ge W \  C_2 e^{\phi(t)} t^{2q+1} \beta(t)  \left(\frac{\epsilon(t)}{\beta(t)}\right)'^2\left(\frac{\beta(t)}{\epsilon(t)}\right)^2.
\end{split} \end{equation}
From \eqref{under_use_1} and \eqref{under_use_2}, it follows that
\[
\frac{d}{dt} \left( e^{\phi(t)} \mathcal{E}(t) \right) \le \frac{C_1}{C_2 W} \frac{d}{dt} \left(e^{\phi(t)} t^{2q+1+m} \beta(t) \left(\frac{\epsilon(t)}{\beta(t)}\right)'^2\left(\frac{\beta(t)}{\epsilon(t)}\right)^2\right).
\]
Integrating the previous differential inequality over $[t_6,t]$, we obtain an explicit upper bound for the Lyapunov functional.
\begin{equation*}
\begin{split}
- e^{\phi(t_6)} \mathcal{E}(t_6) + e^{\phi(t)} \mathcal{E}(t)
&\le \frac{C_1}{C_2 W} \bigg( e^{\phi(t)} t^{2q+1+m} \beta(t)
\left(\frac{\epsilon(t)}{\beta(t)}\right)'^2\left(\frac{\beta(t)}{\epsilon(t)}\right)^2 \\
&\quad - e^{\phi(t_6)} t_6^{2q+1+m} \beta(t_6)
\left(\frac{\epsilon(t_6)}{\beta(t_6)}\right)'^2\left(\frac{\beta(t_6)}{\epsilon(t_6)}\right)^2 \bigg), \\
\mathcal{E}(t)
&\le \frac{C_1}{C_2 W} t^{2q+1+m} \beta(t)
\left(\frac{\epsilon(t)}{\beta(t)}\right)'^2\left(\frac{\beta(t)}{\epsilon(t)}\right)^2 + \frac{C_3}{e^{\phi(t)}}.
\end{split}
\end{equation*}
By the assumption, \( \exists\ \rho_1\in(0,W), \left(\left(\frac{\epsilon(t)}{\beta(t)}\right)'^2\left(\frac{\beta(t)}{\epsilon(t)}\right)^2\right)' \ge -\frac{\rho_1}{t^m} \left(\frac{\epsilon(t)}{\beta(t)}\right)'^2\left(\frac{\beta(t)}{\epsilon(t)}\right)^2,
\)
one obtains for sufficiently large $t$, \(\left(\frac{\epsilon(t)}{\beta(t)}\right)'^2\left(\frac{\beta(t)}{\epsilon(t)}\right)^2 \ge
C e^{-\frac{\rho_1}{1-m} t^{1-m}},
\) consequently, there exists \(C_4 \ge 0\) such that
\[
\mathcal{E}(t) \le C_4 t^{2q+1+m} \beta(t) \left(\frac{\epsilon(t)}{\beta(t)}\right)'^2\left(\frac{\beta(t)}{\epsilon(t)}\right)^2.
\]
We next derive the corresponding estimates for the primal–dual trajectory from this Lyapunov bound. Recalling the definition of the energy functional \eqref{eq:energy}, each component of the energy immediately yields a corresponding estimate for the state variables and velocity terms:
\begin{equation*} \begin{split}
\|x(t)-x_t\|^2 \leq O\left( {t^{2q+m+1}} \beta(t) \left(\frac{\epsilon(t)}{\beta(t)}\right)'^2\left(\frac{\beta(t)}{\epsilon(t)}\right)^2 \right), 
\end{split} \end{equation*}
\begin{equation*} \begin{split}
\|\lambda(t)-\lambda_t\|^2 \leq O\left( {t^{2q+m+1}} \beta(t) \left(\frac{\epsilon(t)}{\beta(t)}\right)'^2\left(\frac{\beta(t)}{\epsilon(t)}\right)^2 \right).
\end{split} \end{equation*}
Moreover, the inertial component can be estimated through the velocity-dependent term contained in $\mathcal{E}(t)$.
\[
\|x(t)-x_t + t^q \dot{x}(t)\|^2 \leq O\left( {t^{2q+m+1}} \beta(t) \left(\frac{\epsilon(t)}{\beta(t)}\right)'^2\left(\frac{\beta(t)}{\epsilon(t)}\right)^2 \right).
\]
In order to recover an explicit decay estimate for the velocity, we combine the inertial component contained in the energy with the previously derived state estimate. Using the elementary inequality
\(
\|u + v\|^2 \le 2\|u\|^2 + 2\|v\|^2,
\)
we further obtain an explicit decay estimate for the velocity $\dot{x}(t)$.
\[
t^{2q} \|\dot{x}(t)\|^2 \le 2\left( \|x(t)-x_t + t^q \dot{x}(t)\|^2 + \|x(t)-x_t\|^2 \right).
\]
Substituting the existing estimates
\(
\|x(t)-x_t + t^q \dot{x}(t)\|^2 \le 2C_4 t^{2q+1+m} \beta(t)\left(\frac{\epsilon(t)}{\beta(t)}\right)'^2\left(\frac{\beta(t)}{\epsilon(t)}\right)^2\
\)
and
\( \|x(t)-x_t\|^2 \le C_4 t^{2q+1+m} \beta(t) \left(\frac{\epsilon(t)}{\beta(t)}\right)'^2\left(\frac{\beta(t)}{\epsilon(t)}\right)^2,
\)
we obtain
\[
t^{2q} \|\dot{x}(t)\|^2 \le 6 C_4 t^{2q+1+m} \beta(t)\left(\frac{\epsilon(t)}{\beta(t)}\right)'^2\left(\frac{\beta(t)}{\epsilon(t)}\right)^2 ,
\]
so
\[
\|\dot{x}(t)\|^2 \le O\left( t^{m+1} \beta(t) \left(\frac{\epsilon(t)}{\beta(t)}\right)'^2\left(\frac{\beta(t)}{\epsilon(t)}\right)^2\right) .
\]
We next estimate the feasibility violation, applying Lemma \ref{lem2.3}(i) yields
\[
\|Ax(t)-b\| \leq O\left( \frac{\epsilon(t)}{\beta(t)} + {t^{\frac{2q+m+1}{2}}} \sqrt{\beta(t)}\left|\left(\frac{\epsilon(t)}{\beta(t)}\right)'\right| \frac{\beta(t)}{\epsilon(t)} \right).
\]
Also
\[
\|A(x(t)+\Gamma(t) \dot{x}(t))-b\| \le \|Ax(t)-b\| + \|A\| \|\Gamma(t)|\|\dot{x}(t)\|.
\]
We now split the discussion into two cases depending on the values $r$ and $r_0$. If \(r_0 = 0, \ r \ge 0\), \(|\Gamma(t)|  \|\dot{x}(t)\| \le O\left( {t^{\frac{m+1-2q}{2}}} \sqrt{\beta(t)} \left|\left(\frac{\epsilon(t)}{\beta(t)}\right)'\right| \frac{\beta(t)}{\epsilon(t)}\right),\) 
If \(r_0 > 0, \ r \in \mathbb{R}\), \(|\Gamma(t)|  \|\dot{x}(t)\|  \le O\left( {t^{\frac{m+1}{2}}} \sqrt{\beta(t)} \left|\left(\frac{\epsilon(t)}{\beta(t)}\right)'\right| \frac{\beta(t)}{\epsilon(t)}\right).\) 
As a consequence, \(|\Gamma(t)|  \|\dot{x}(t)\| \le O\left( {t^{\frac{m+1}{2}}} \sqrt{\beta(t)} \left|\left(\frac{\epsilon(t)}{\beta(t)}\right)'\right| \frac{\beta(t)}{\epsilon(t)}\right)\). Since \(  {t^{\frac{m+1-2q}{2}}}  < {t^{\frac{m+1}{2}}} \),
hence the contribution induced by the inertial correction term is of lower asymptotic order and can be absorbed into the leading term. Hence, \[\|A(x(t)+\Gamma(t) \dot{x}(t))-b\| \le O\left( \frac{\epsilon(t)}{\beta(t)} +  {t^{\frac{2q+m+1}{2}}} \sqrt{\beta(t)} \left|\left(\frac{\epsilon(t)}{\beta(t)}\right)'\right| \frac{\beta(t)}{\epsilon(t)}\right).\]
Moreover, by the definition of Lyapunov function\eqref{eq:energy}, it yields
\[
\mathcal{L}_{\epsilon(t)}(x(t)+\Gamma(t)\dot{x}(t), \lambda_t) - \mathcal{L}_{\epsilon(t)}(x_t, \lambda_t) \le O\left({t^{m+1}}\left(\frac{\epsilon(t)}{\beta(t)}\right)'^2\left(\frac{\beta(t)}{\epsilon(t)}\right)^2 \right),
\]

Combining Lemma \ref{lem2.3}, we obtain
\begin{equation} \begin{split}\label{L-L}
\mathcal{L}\bigl(x(t) + \Gamma(t)\dot{x}(t), \lambda^*\bigr) - \mathcal{L}(x^*, \lambda^*)
\le \mathcal{L}_{\epsilon(t)}\bigl(x(t) + \Gamma(t)\dot{x}(t), \lambda_t\bigr) - \mathcal{L}_{\epsilon(t)}(x_t, \lambda_t) \\
\quad + \|\lambda_t - \lambda^*\| \left\| A\bigl(x(t) + \Gamma(t)\dot{x}(t)\bigr) - b \right\|
+ \frac{\epsilon(t)}{2\beta(t)} \ \|x^*\|^2  ,
\end{split}\end{equation} 
further we can get
\begin{equation*}\begin{split} 
&\mathcal{L}(x(t)+\Gamma(t)\dot{x}(t), \lambda^*) - \mathcal{L}(x^*, \lambda^*)\\
&\leq O\left( \frac{\epsilon(t)}{\beta(t)} + {t^{\frac{2q+m+1}{2}}} \left|\left(\frac{\epsilon(t)}{\beta(t)}\right)'\right| \frac{\beta(t)}{\epsilon(t)}\sqrt{\beta(t)}  +{t^{m+1}}\left(\frac{\epsilon(t)}{\beta(t)}\right)'^2\left(\frac{\beta(t)}{\epsilon(t)}\right)^2 \right).
\end{split} \end{equation*}
Finally, using the definition of the Lagrangian function together with the feasibility estimate derived above, we obtain the convergence rate for the objective residual,
\[
|f(x(t)+\Gamma(t) \dot{x}(t)) - f(x^*)  | \le L(x(t)+\Gamma(t) \dot{x}(t), \lambda^*) - L(x^*, \lambda^*) + \|\lambda^*\| \ \|A(x(t)+\Gamma(t) \dot{x}(t))-b\|.
\]
Thus, combining the Lagrangian identity with the feasibility estimate derived above, we obtain 
\begin{equation*}\begin{split}
&|f(x(t)+\Gamma(t) \dot{x}(t)) - f(x^*)| \\
&\leq O\left( \frac{\epsilon(t)}{\beta(t)} + {t^{\frac{2q+m+1}{2}}} \left|\left(\frac{\epsilon(t)}{\beta(t)}\right)'\right| \frac{\beta(t)}{\epsilon(t)} \sqrt{\beta(t)} +{t^{m+1}}\left(\frac{\epsilon(t)}{\beta(t)}\right)'^2\left(\frac{\beta(t)}{\epsilon(t)}\right)^2 \right).
\end{split} \end{equation*}

\paragraph{Case (II):}
\[
t^{2q} \beta(t)\left| \left(\frac{\epsilon(t)}{2\beta(t)}\right)'\right| \ge \max\left\{ t^{2q+1} \beta(t)\left(\frac{\epsilon(t)}{\beta(t)}\right)'^2\left(\frac{\beta(t)}{\epsilon(t)}\right)^2, \frac{1}{t^{q} \epsilon(t)}\left(\frac{\epsilon(t)}{\beta(t)}\right)'^2\left(\frac{\beta(t)}{\epsilon(t)}\right)^2 \right\}
\]
Repeating the same comparison strategy as in previous cases, we obtain
\[
\frac{d}{dt} \left( e^{\phi(t)} \mathcal{E}(t) \right) \le - \widehat{C}_1
e^{\phi(t)}t^{2q} \beta(t) \left(\frac{\epsilon(t)}{\beta(t)}\right)'.
\]
Let $t_7:=\max\{T_0,T_2\},\ \widehat{C}_2 :=\frac{W-\rho_2}{W}\in(0,1).$
Direct computation yields
\begin{equation*}
\begin{split}
\frac{d}{dt} \left( e^{\phi(t)} t^{2q} \beta(t) t^m \ \left(-\frac{\epsilon(t)}{\beta(t)}\right)' \right)
&= \frac{d}{dt} \left( e^{\phi(t)} t^{2q+m} \beta(t) \ \left(-\frac{\epsilon(t)}{\beta(t)}\right)' \right) \\
&= e^{\phi(t)} \left( \frac{W}{t^m} \ t^{2q+m}\beta(t)  \right) \left(-\frac{\epsilon(t)}{\beta(t)}\right)' \\
&\quad+ e^{\phi(t)}(2q+m) t^{2q+m-1}\beta(t) \left(-\frac{\epsilon(t)}{\beta(t)}\right)' \\
&\quad + e^{\phi(t)}t^{2q+m}{\beta(t)}' \left(-\frac{\epsilon(t)}{\beta(t)}\right)'+e^{\phi(t)}t^{2q+m}\beta(t)\ \left(-\frac{\epsilon(t)}{\beta(t)}\right)''\\
&\ge W \widehat{C}_2 t^{2q}\beta(t)e^{\phi(t)} \left(-\frac{\epsilon(t)}{\beta(t)}\right)' .
\end{split}
\end{equation*}
Combining the above two inequalities, we deduce that
\[
\frac{d}{dt} \left( e^{\phi(t)} \mathcal{E}(t) \right) \le \frac{\widehat{C}_1}{\widehat{C}_2 W} \frac{d}{dt} \left(e^{\phi(t)} t^{2q+m} \beta(t) \left(-\frac{\epsilon(t)}{\beta(t)}\right)'\right).
\]
Integrating both sides for \(t \ge t_7\), we obtain
\begin{equation*} \begin{split}
- e^{\phi(t_7)} \mathcal{E}(t_7) + e^{\phi(t)} \mathcal{E}(t)
&\le \frac{\widehat{C}_1}{\widehat{C}_2 W} \left( e^{\phi(t)} t^{2q+m} \beta(t) \left(-\frac{\epsilon(t_7)}{\beta(t)}\right)' - e^{\phi(t_7)} t_7^{2q+m} \beta(t) \left(-\frac{\epsilon(t_7)}{\beta(t_7)}\right)' \right), \\
 \mathcal{E}(t)
&\le \frac{\widehat{C}_1}{\widehat{C}_2 W} t^{2q+m} \beta(t) \left(-\frac{\epsilon(t)}{\beta(t)}\right)' + \frac{\widehat{C}_3}{e^{\phi(t)}}.
\end{split} \end{equation*}
As in the previous case, there exists \(\widehat{C}_4 \ge 0\) such that
\[
\mathcal{E}(t) \le \widehat{C}_4 t^{2q+m} \beta(t) \left(-\frac{\epsilon(t)}{\beta(t)}\right)'.
\]
Similarly, we obtain
\begin{equation*} \begin{split}
\|x(t)-x_t\|^2 &\leq O\left( t^{2q+m} \beta(t) \left(-\frac{\epsilon(t)}{\beta(t)}\right)' \right) ,\quad
\|\lambda(t)-\lambda_t\|^2 \leq  O\left( t^{2q+m} \beta(t) \left(-\frac{\epsilon(t)}{\beta(t)}\right)' \right)  
\end{split} \end{equation*}
\begin{equation*} \begin{split}
t^{2q} \|\dot{x}(t)\|^2 &\le 6 \widehat{C}_4 t^{2q+m} \beta(t) \left(-\frac{\epsilon(t)}{\beta(t)}\right)', 
\end{split} \end{equation*}
so we can get
\begin{equation*} \begin{split}
\|\dot{x}(t)\|^2 &\leq O\left( t^{m} \beta(t) \left(-\frac{\epsilon(t)}{\beta(t)}\right)' \right) .
\end{split} \end{equation*}
Applying Lemma \ref{lem2.3}(i), we obtain
\[
\|Ax(t)-b\| \leq O\left( \frac{\epsilon(t)}{\beta(t)} + {t^{\frac{2q+m}{2}}} \sqrt{ \beta(t)\left(-\frac{\epsilon(t)}{\beta(t)}\right)'} \right).
\]
Moreover,
\[
\|A(x(t)+\Gamma(t) \dot{x}(t))-b\| \le \|Ax(t)-b\| + \|A\| \ |\Gamma(t)| \ \|\dot{x}(t)\|.
\]
We consider two separate cases for $r$ and $r_0$. \\
If \(r_0 = 0, r \ge 0\), \(|\Gamma(t)| \ \|\dot{x}(t)\| \le O\left( {t^{\frac{m-2q}{2}}} \sqrt { \beta(t)\left(-\frac{\epsilon(t)}{\beta(t)}\right)'} \right)\),
If \(r_0 > 0, r \in \mathbb{R}\), \(|\Gamma(t)| \ \|\dot{x}(t)\| \le O\left({t^{\frac{m}{2}}} \sqrt{\beta(t)\left( -\frac{\epsilon(t)}{\beta(t)}\right)'} \right)\),
In summary, \(|\Gamma(t)| \ \|\dot{x}(t)\| \le O\left({t^{\frac{m}{2}}} \sqrt{ \beta(t)\left(-\frac{\epsilon(t)}{\beta(t)}\right)'} \right)\). Since \({t^{\frac{2q+m}{2}}} \) > \({t^{\frac{m}{2}}} \) , the former term dominates asymptotically. Hence,
\[
 \|A(x(t)+\Gamma(t) \dot{x}(t))-b\| \leq O\left( \frac{\epsilon(t)}{\beta(t)} + {t^{\frac{2q+m}{2}}} \sqrt{\beta(t)\left(-\frac{\epsilon(t)}{\beta(t)}\right)'} \right).
\]
Furthermore,
\begin{equation*} \begin{split}
\mathcal{L}_{\epsilon(t)}(x(t)+\Gamma(t)\dot{x}(t), \lambda_t) - \mathcal{L}_{\epsilon(t)}(x_t, \lambda_t)
& \leq O\left( {t^{2q+m}}\beta(t) \left(-\frac{\epsilon(t)}{\beta(t)}\right)' \ \frac{1}{t^{2q} \beta(t)} \right) \\
&= O\left({t^{m}} \left(-\frac{\epsilon(t)}{\beta(t)}\right)'  \right).
\end{split} \end{equation*}

By the same reasoning, we can derive
\[
\mathcal{L}(x(t)+\Gamma(t)\dot{x}(t), \lambda^*) - \mathcal{L}(x^*, \lambda^*) \leq O\left( \frac{\epsilon(t)}{\beta(t)} + {t^{\frac{2q+m}{2}}} \sqrt{\beta(t)\left(-\frac{\epsilon(t)}{\beta(t)}\right)'} + {t^{m}} \left(-\frac{\epsilon(t)}{\beta(t)}\right)' \right).
\]
By the definition of \(\mathcal{L}\),
\[
|f(x(t)+\Gamma(t) \dot{x}(t)) - f(x^*)| \le \mathcal{L}(x(t)+\Gamma(t) \dot{x}(t), \lambda^*) - \mathcal{L}(x^*, \lambda^*) + \|\lambda^*\| \ \|A(x(t)+\Gamma(t) \dot{x}(t))-b\|,
\]
so
\[
|f(x(t)+\Gamma(t) \dot{x}(t)) - f(x^*)| \leq O\left( \frac{\epsilon(t)}{\beta(t)} + {t^{\frac{2q+m}{2}}} \sqrt{\beta(t)\left(-\frac{\epsilon(t)}{\beta(t)}\right)'} + {t^{m} } \left(-\frac{\epsilon(t)}{\beta(t)}\right)' \right).
\]
In both cases, the preceding estimates together with the corresponding limit assumptions imply
\[
\lim_{t \to +\infty} \|(x(t), \lambda(t)) - (x_t, \lambda_t)\| = \lim_{t \to +\infty}\sqrt{\|x(t)-x_t\|^2 + \|\lambda(t)-\lambda_t\|^2} = 0.
\]
On the other hand, Lemma \ref{lem2.1} gives
\[
\lim_{t \to +\infty} \|(x_t, \lambda_t) - (x^*, \lambda^*)\| = 0.
\]
Therefore, by the triangle inequality,
\[
\|(x(t), \lambda(t)) - (x^*, \lambda^*)\| \le \|(x(t), \lambda(t)) - (x_t, \lambda_t)\| + \|(x_t, \lambda_t) - (x^*, \lambda^*)\|,
\]
and hence
\[
\lim_{t \to +\infty} \|(x(t), \lambda(t)) - (x^*, \lambda^*)\| = 0.
\]
Thus, the whole primal‑dual trajectory converges strongly to the minimum‑norm element of the saddle‑point set. This completes the proof.
\end{proof}

\textbf{Remark 3.2} The two sets of assumptions in Theorem 3.2 are nonempty and can be satisfied by commonly used time-dependent coefficient functions. 

For Case (i), consider the power-type coefficients$\beta(t)=t^{s},\ \epsilon(t)=t^{-p},$ where $s,p>0$. If $0<q<m<1,\ q+s+p<1,\ p<q+m,\ 2q+m+s<1,$ then the assumptions in \eqref{H0} are satisfied. Moreover,
\[
\left(\frac{\epsilon(t)}{\beta(t)}\right)'^2 \left(\frac{\beta(t)}{\epsilon(t)}\right)^2=\frac{(p+s)^2}{t^2},
\]
and hence, for every $\rho_1\in(0,W)$, the first differential condition in Case (i) holds for all sufficiently large $t$. Furthermore,
\[\frac{t^{2q+1}\beta(t)\left(\frac{\epsilon(t)}{\beta(t)}\right)'^2\left(\frac{\beta(t)}{\epsilon(t)}\right)^2}{t^{2q}\beta(t)\left|\left(\frac{\epsilon(t)}{2\beta(t)}\right)'\right|}=2(p+s)t^{p+s}\longrightarrow +\infty,\]
while
\[t^{2q+1+m}\beta(t)\left(\frac{\epsilon(t)}{\beta(t)}\right)'^2\left(\frac{\beta(t)}{\epsilon(t)}\right)^2=(p+s)^2 t^{2q+m+s-1}\longrightarrow 0.\]
Thus, all the assumptions of Case (i) are fulfilled.

For Case (ii), non-power-type coefficients can also be admitted. Let
$\beta(t)=\big(\log(t+e)\big)^s,\ \epsilon(t)=\big(\log(t+e)\big)^{-p},$
where $s,p>0$, and set $a=p+s$. Suppose that
$0<a<1,\ 0<q<m<1,\ 2q+m<1.$
Then
\[\frac{\epsilon(t)}{\beta(t)}=\big(\log(t+e)\big)^{-a}, \quad 
\left(-\frac{\epsilon(t)}{\beta(t)}\right)'=\frac{a}{(t+e)\big(\log(t+e)\big)^{a+1}}.\]
Since $m<1$,
\[
\frac{\left(-\frac{\epsilon(t)}{\beta(t)}\right)''}{\left(-\frac{\epsilon(t)}{\beta(t)}\right)'}=-\frac{1}{t+e}-\frac{a+1}{(t+e)\log(t+e)}=o(t^{-m}),
\]
and therefore, for every $\rho_2\in(0,W)$,
$\left(-\frac{\epsilon(t)}{\beta(t)}\right)''\ge -\frac{\rho_2}{t^m}\left(-\frac{\epsilon(t)}{\beta(t)}\right)'$
for all sufficiently large $t$. In addition,
\[
\frac{t^{2q}\beta(t)\left|\left(\frac{\epsilon(t)}{2\beta(t)}\right)'\right|}{t^{2q+1}\beta(t)\left(\frac{\epsilon(t)}{\beta(t)}\right)'^2\left(\frac{\beta(t)}{\epsilon(t)}\right)^2}\sim \frac{(\log t)^{1-a}}{2a}\longrightarrow +\infty, 
\]
\[
t^{2q+m}\beta(t)\left|\left(\frac{\epsilon(t)}{2\beta(t)}\right)'\right|=O\big(t^{2q+m-1}(\log t)^{-p-1}\big)\longrightarrow 0.
\]
Hence the assumptions of Case (ii) are also realizable.

\textbf{Remark 3.3} In the present dynamical system, the time scaling coefficient functions \(\beta(t)\) and Tikhonov regularization coefficient \(\epsilon(t)\) are not restricted to prescribed power-law forms. This general family of time‑dependent coefficients encompasses functions with qualitatively distinct asymptotic behaviors. Thus, the main feature of Theorem \ref{Th3.2} is the simultaneous realization of a mixed-order structure, implicit Hessian damping, general coefficient functions, and minimum-norm primal-dual selection within a single dynamical framework.

\section{Numerical experiments}\label{sec:num}

Numerical experiments are presented to illustrate the convergence behavior of the proposed dynamical system, the oscillation suppression effect of the implicit Hessian damping mechanism, the minimum-norm selection induced by Tikhonov regularization, and the numerical performance of the proposed method in comparison with existing primal-dual dynamical systems. All experiments were performed on a Lenovo ThinkBook 16 G5+ IRH (16GB) platform using MATLAB R2025a, and the dynamical systems were numerically solved with the built-in solver ode23. 

Example 4.1.  \ We first consider a convex but not strictly convex constrained optimization problem:
\[
\min_{x\in\mathbb{R}^3} f(x)=\frac{1}{2}(x_1^2+(x_2-x_3)^2)\quad\text{s.t.}\quad 2x_1-x_2+x_3=2.
\]
The solution set of the constrained optimization problem is
\[
S=\left\{(x_1,x_2,x_3)\in\mathbb{R}^3\;\middle|\;x_1=\frac{4}{5},\;x_2-x_3=-\frac{2}{5}\right\},
\]
\noindent and the minimum-norm solution is \(x^*=(\frac{4}{5},-\frac{1}{5},\frac{1}{5})^\top\).These parameter choices are compatible with the assumptions of Theorems \ref{Th3.1} and \ref{Th3.2}, for instance, one may take $m=0.25$ and $W=0.2$. In this first experiment, we set \(q=0.2\) , \(\alpha=6\), and consider four groups of comparative parameters:
\[
\{\beta(t),\epsilon(t),r_0,r\}=
\left\{
(t^{0.2}, t^{-0.3}, 0, 0.3);\;
(t^{0.2}, t^{-0.4}, 0, 0.3);\;
(t^{0.3}, t^{-0.3}, 0, 0.3);\;
(t^{0.2}, t^{-0.3}, 0, 0)
\right\}.
\]
The initial conditions are \(x(1)=(1,1,1)^\top,\;\lambda(1)=1,\;\dot{x}(1)=(1,1,1)^\top\). Particularly the last set of parameters corresponds to the dynamical system without the implicit Hessian term.

\begin{figure}[htbp]
\centering
\begin{tabular}{ccc}
\begin{subfigure}{0.31\textwidth}
    \centering
    \includegraphics[width=\linewidth]{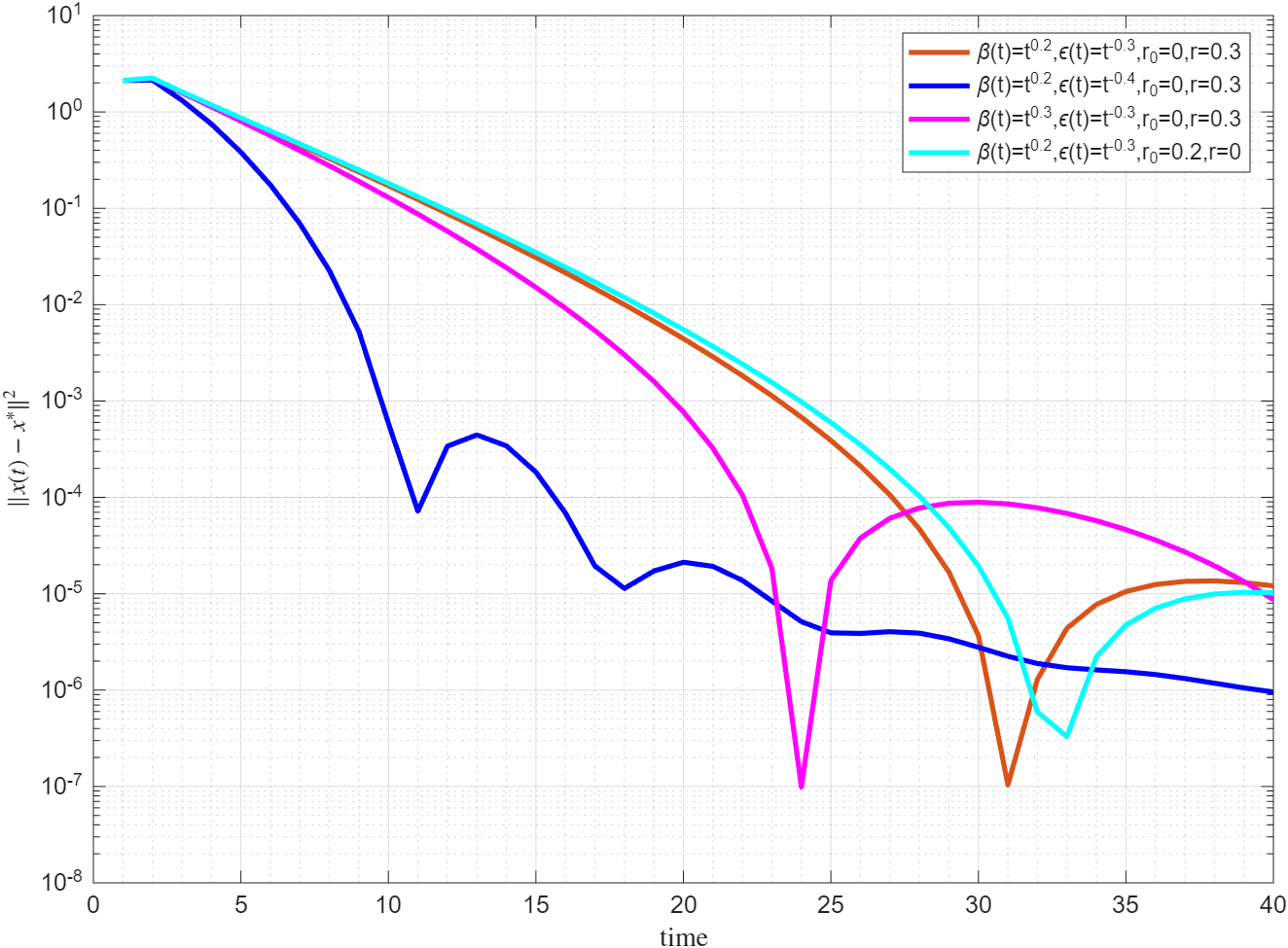}
    \caption{iterate error}
    \label{my:x}
\end{subfigure}
&
\begin{subfigure}{0.31\textwidth}
    \centering
    \includegraphics[width=\linewidth]{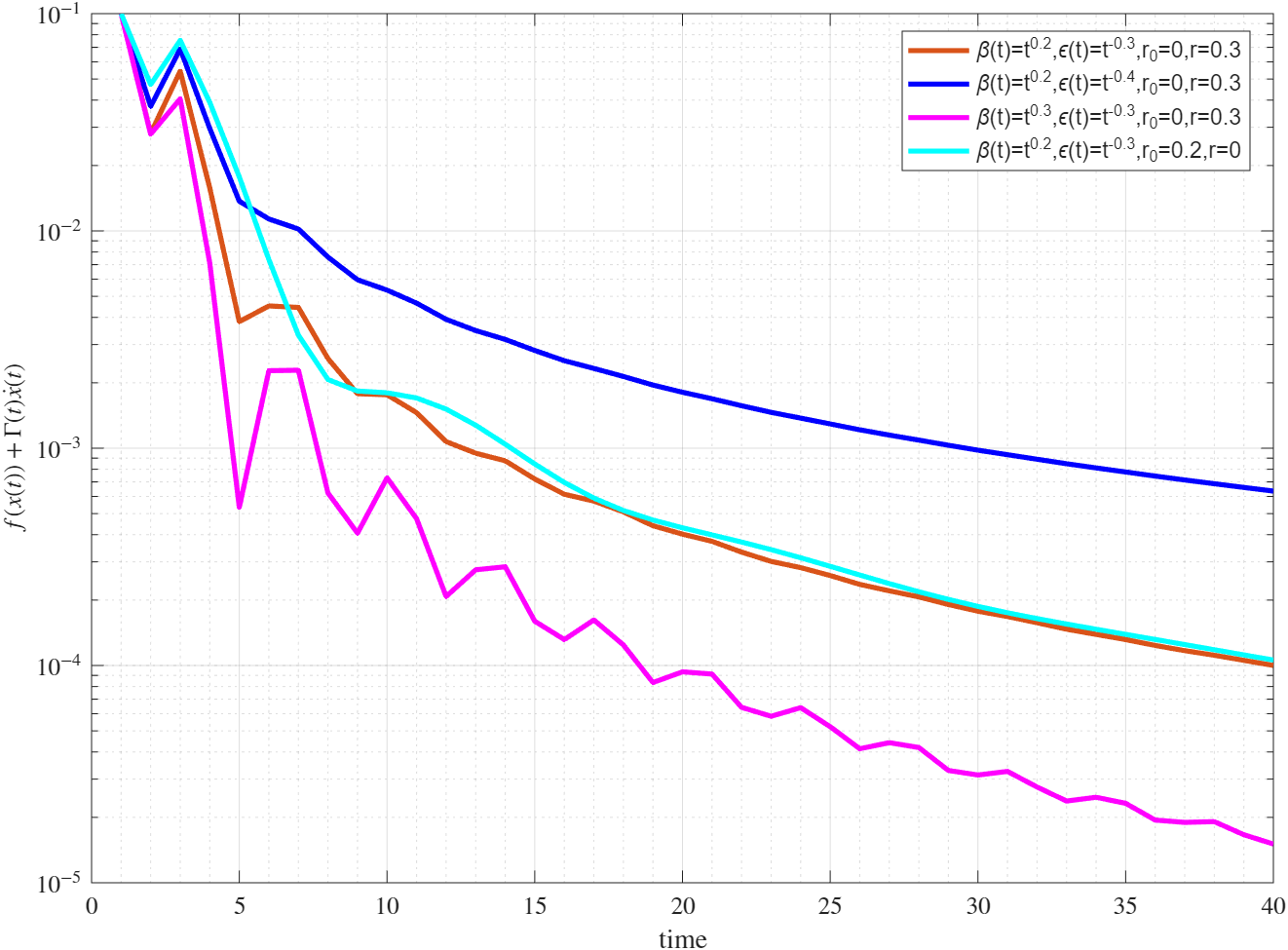}
    \caption{objective gap}
    \label{my:f}
\end{subfigure}
&
\begin{subfigure}{0.31\textwidth}
    \centering
    \includegraphics[width=\linewidth]{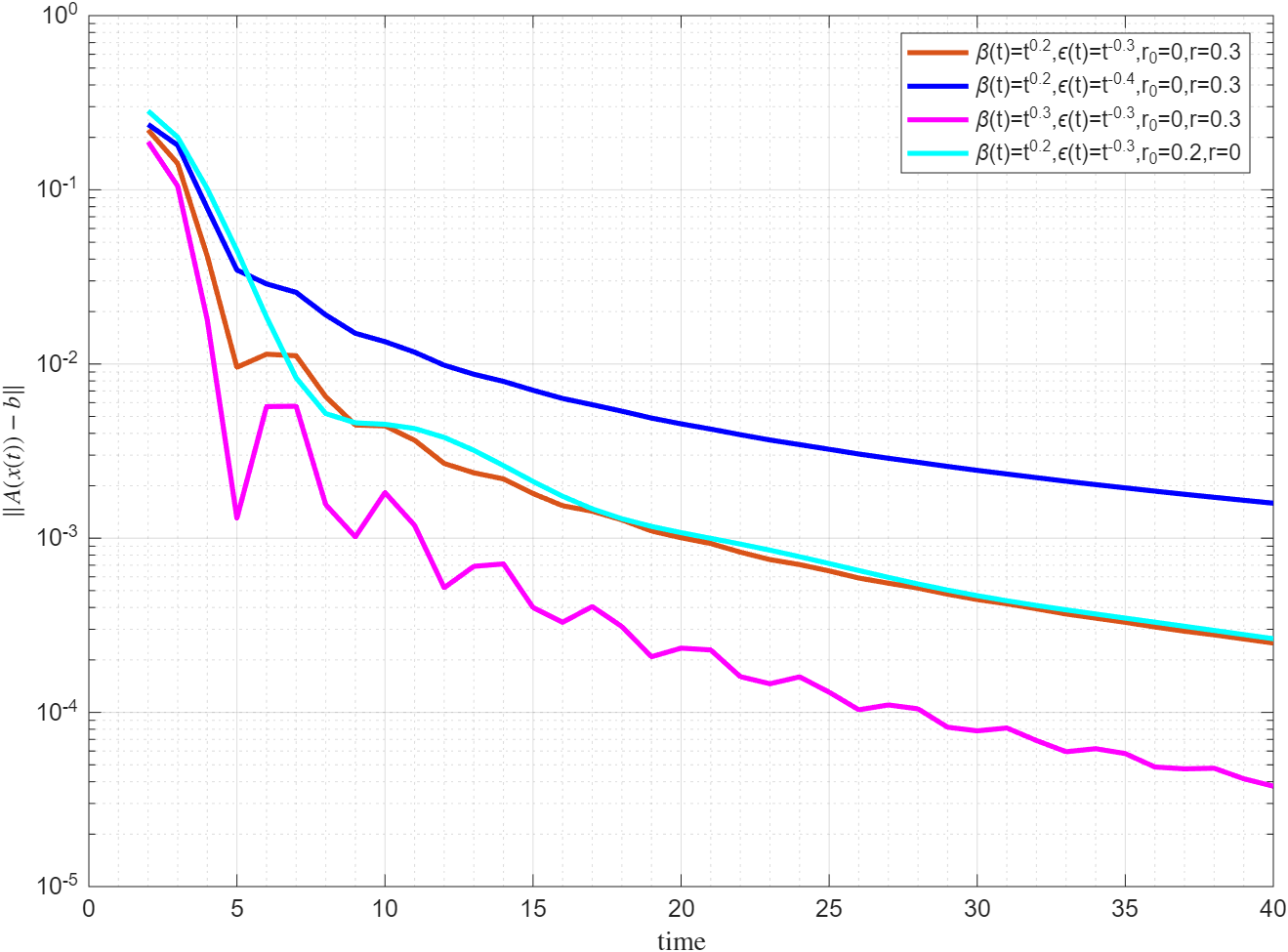}
    \caption{feasibility violation}
    \label{my:Ax-b}
\end{subfigure}
\end{tabular}

\caption{Errors of  the trajectory, the objective gap, and the feasibility violation of our dynamical system}
\label{fig:F1}
\end{figure}

Figure \ref{fig:F1} compares the convergence performance of \(\|x(t)-x^*\|^2\), \(|f(x(t))-f(x^*)|\) and \(\|Ax(t)-b\|\) of the proposed primal-dual dynamical system  for solving the constrained optimization problem under four different parameter settings. 

\begin{figure}[H]
\centering
\begin{tabular}{cc}
\begin{subfigure}{0.45\textwidth}
    \centering
    \includegraphics[width=\linewidth]{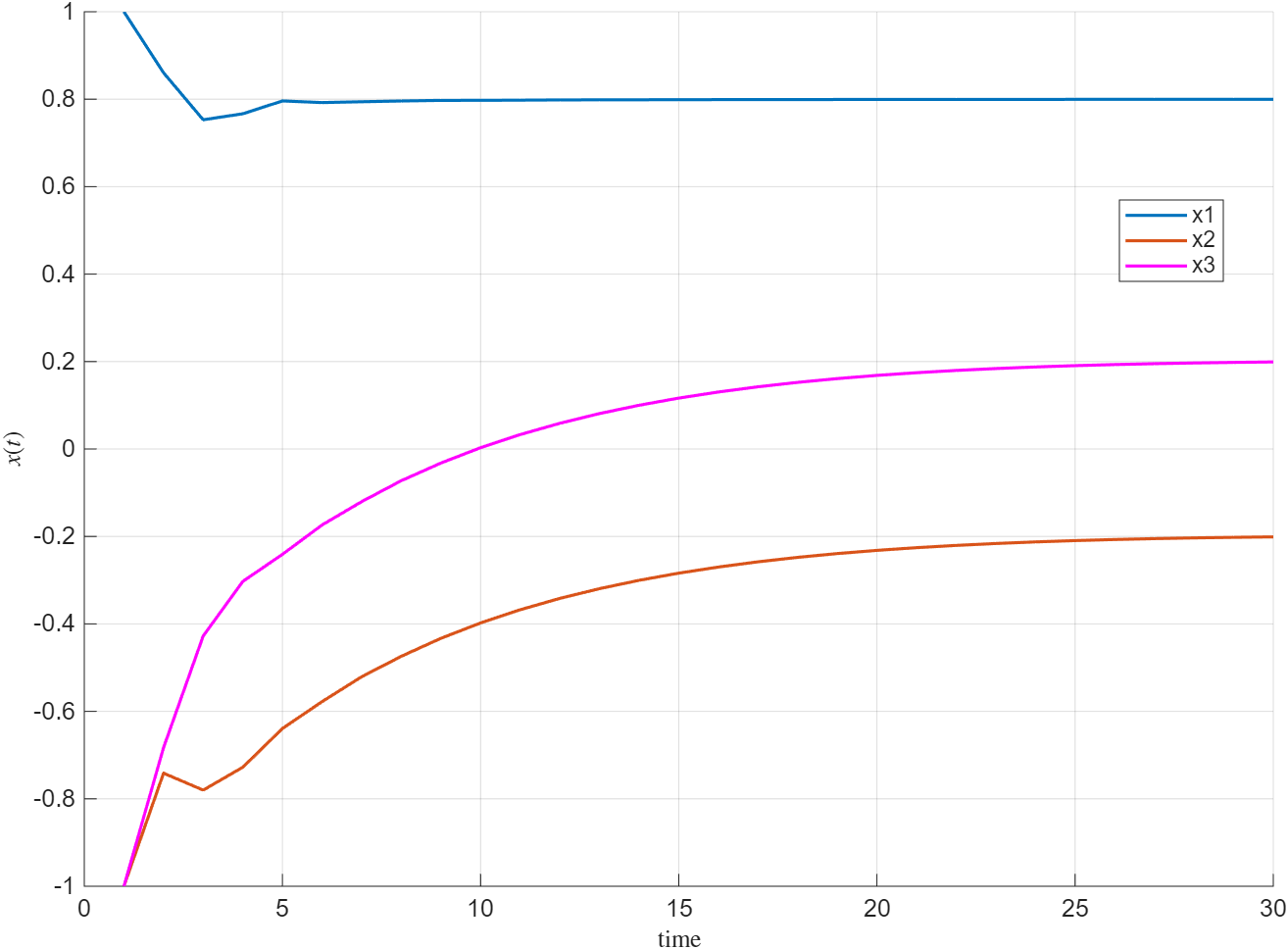}
    \caption{trajectory with Tikhonov regularization}
    \label{withTik}
\end{subfigure}
&
\begin{subfigure}{0.45\textwidth}
    \centering
    \includegraphics[width=\linewidth]{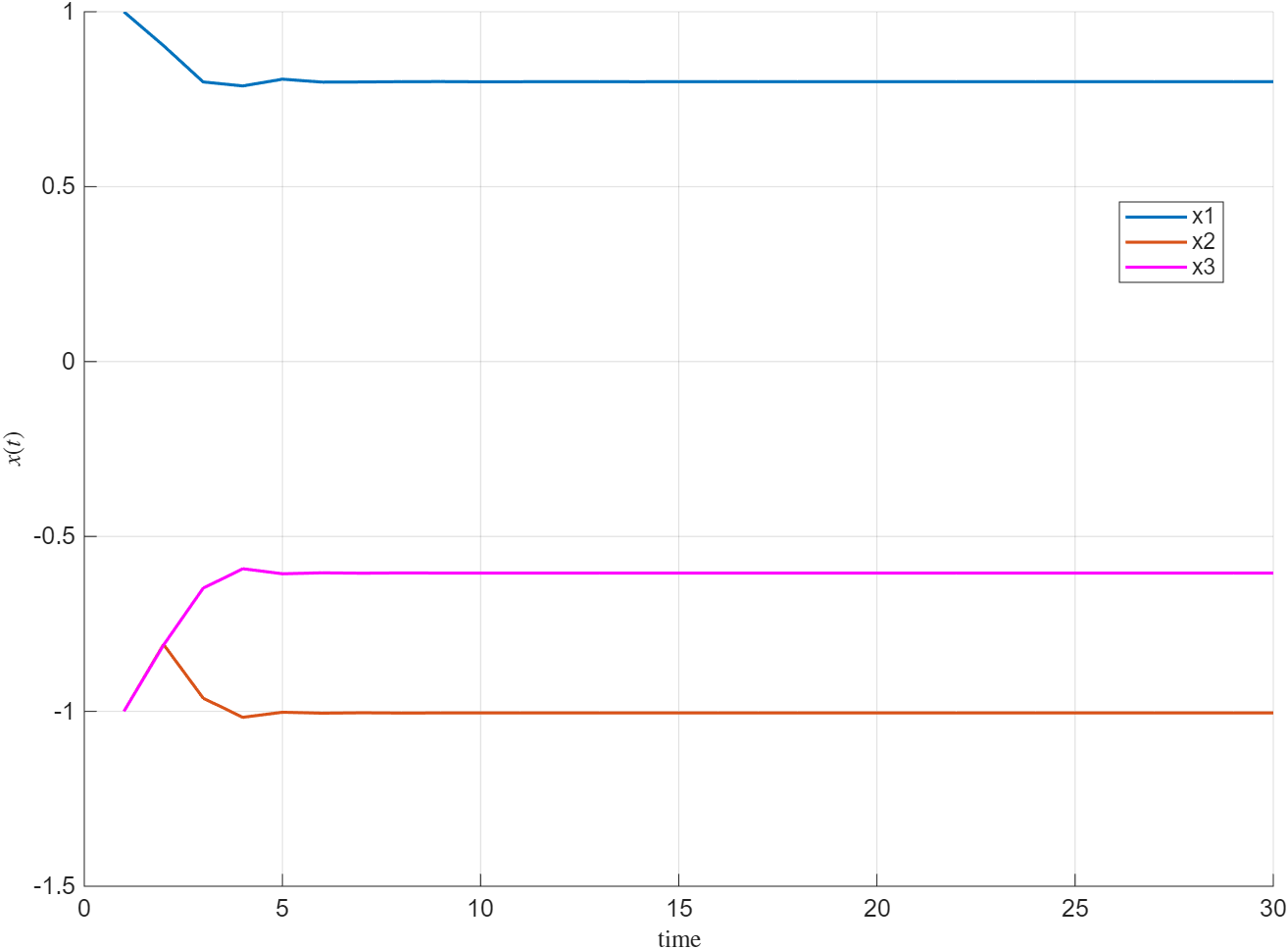}
    \caption{trajectory without Tikhonov regularization}
    \label{withoutTik}
\end{subfigure}
\end{tabular}
\caption{Comparison of trajectory components of the dynamical system with and without Tikhonov regularization}
\label{fig:F2}
\end{figure}

To highlight the role of Tikhonov regularization, we also compare the dynamical system without the Tikhonov regularization term (\(\epsilon(t)\equiv0\)). As displayed in Figure \ref{fig:F2}, the trajectory with the Tikhonov regularization term converges to the minimum-norm solution \(x^*\), whereas the trajectory without Tikhonov regularization approaches a different solution and does not exhibit the minimum-norm selection property.

\begin{figure}[H]
\centering
\begin{tabular}{ccc}
\begin{subfigure}{0.3\textwidth}
    \centering
    \includegraphics[width=\linewidth, keepaspectratio]{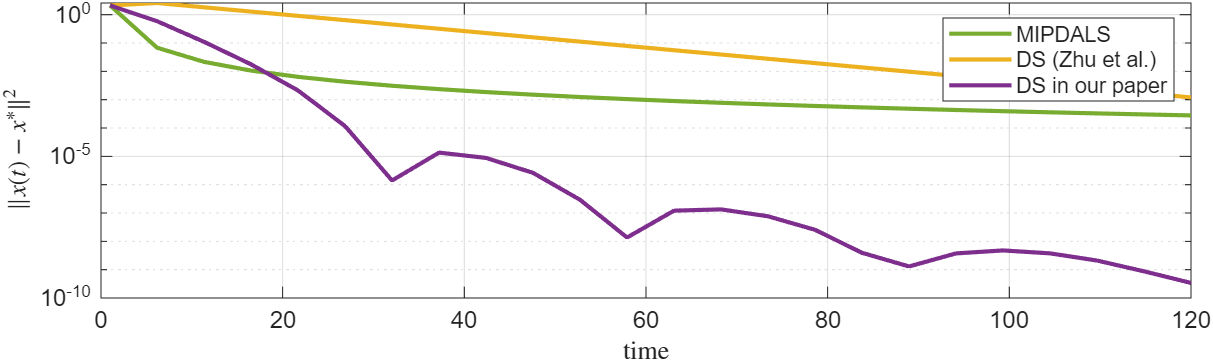}
    \caption{iterate error}
    \label{com:x}
\end{subfigure}
&
\begin{subfigure}{0.3\textwidth}
    \centering
    \includegraphics[width=\linewidth, keepaspectratio]{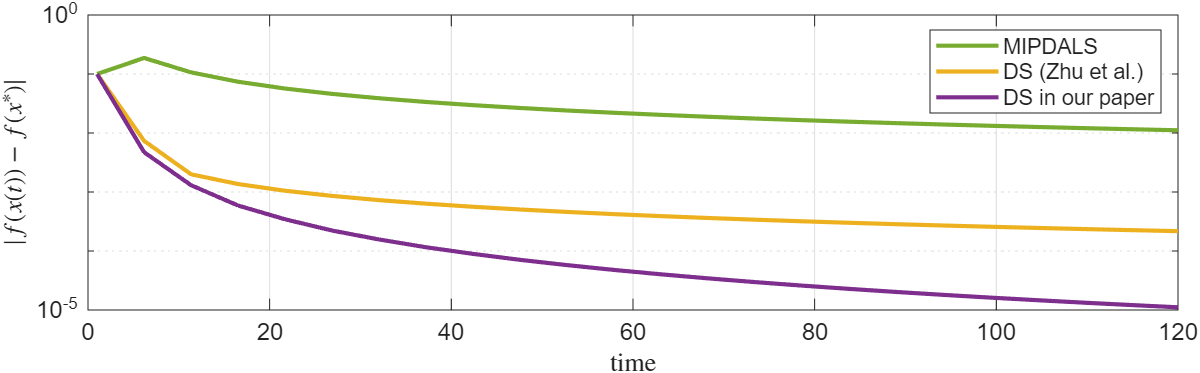}
    \caption{objective gap}
    \label{com:f}
\end{subfigure}
&
\begin{subfigure}{0.3\textwidth}
    \centering
    \includegraphics[width=\linewidth, keepaspectratio]{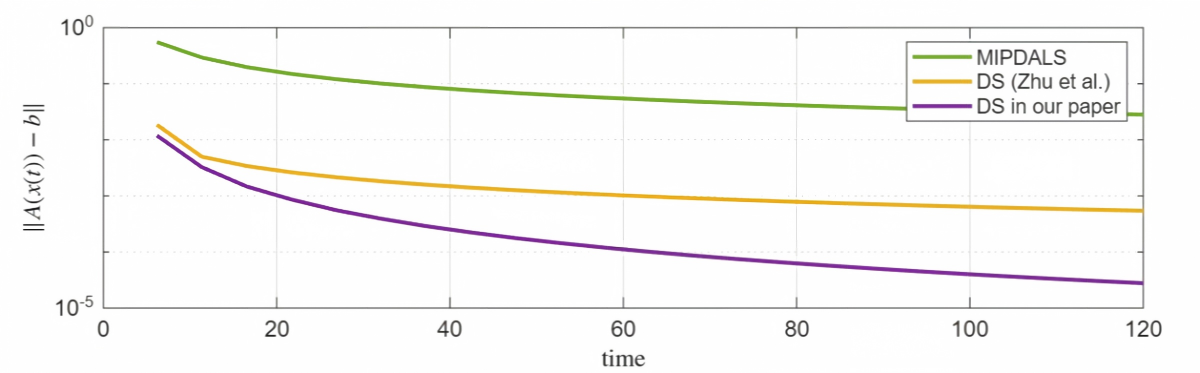}
    \caption{feasibility violation}
    \label{com:Ax-b}
\end{subfigure}
\end{tabular}
\caption{Comparison of error results between MIPDALS, the dynamical system proposed by Zhu et al. and our dynamical system}
\label{fig:F3}
\end{figure}

In the second experiment, we compare our dynamical system with two existing dynamical systems: the MIPDALS dynamical system (Battahi et al. \cite{fouad_battahi_simultaneous_2025}) and the dynamical system proposed by Zhu et al \cite{zhu_strong_2024}. The parameters of the comparison dynamical systems were selected according to the numerical settings used in the corresponding references. Figure \ref{fig:F3} compares the numerical performance of the three dynamical systems in terms of the iterate error, objective residual evaluated along the primal trajectory, and feasibility violation. Under the parameter settings used in the experiment, the proposed dynamical system exhibits faster decay in all three metrics.

Example 4.2. \ We consider an additional convex optimization problem:

\begin{equation*}
\min_{x\in\mathbb{R}^3} \frac14\bigl(x_1+2x_2-x_3\bigr)^4+\frac{1}{16}\bigl(2x_1-x_2+x_3-3\bigr)^4+\log\cosh\bigl(4(x_1-x_2-2x_3+1)\bigr)\quad \text{s.t.}\quad x_1+x_2-x_3=0,
\end{equation*}
An easy computation shows that the set of solutions for this problem is given by
\(
S=\Bigl\{\bigl(1,0,1)^T\Bigr\},
\)
and the minimum-norm solution is $x^*=(1,0,1)^T$. We perform this experiment in which all parameters are fixed except for $\Gamma(t)=r_0+\frac{r}{t^q}$. 
We take
\(
\beta(t)=t^{0.2},\ \epsilon(t)=t^{-0.3},
\)
and compare the following three parameter settings:
\(
(r_0,r)=\{(0,2),(0,0),(2,0)\}
\)
The choice $(r_0,r)=(0,0)$ gives $\Gamma(t)\equiv0$ and therefore corresponds to the dynamical 
system without implicit Hessian damping, whereas the other two choices yield nonzero extrapolation 
terms and activate the implicit Hessian damping.

\begin{figure}[htbp]
\centering
\begin{tabular}{ccc}
\begin{subfigure}{0.31\textwidth}
    \centering
    \includegraphics[width=\linewidth]{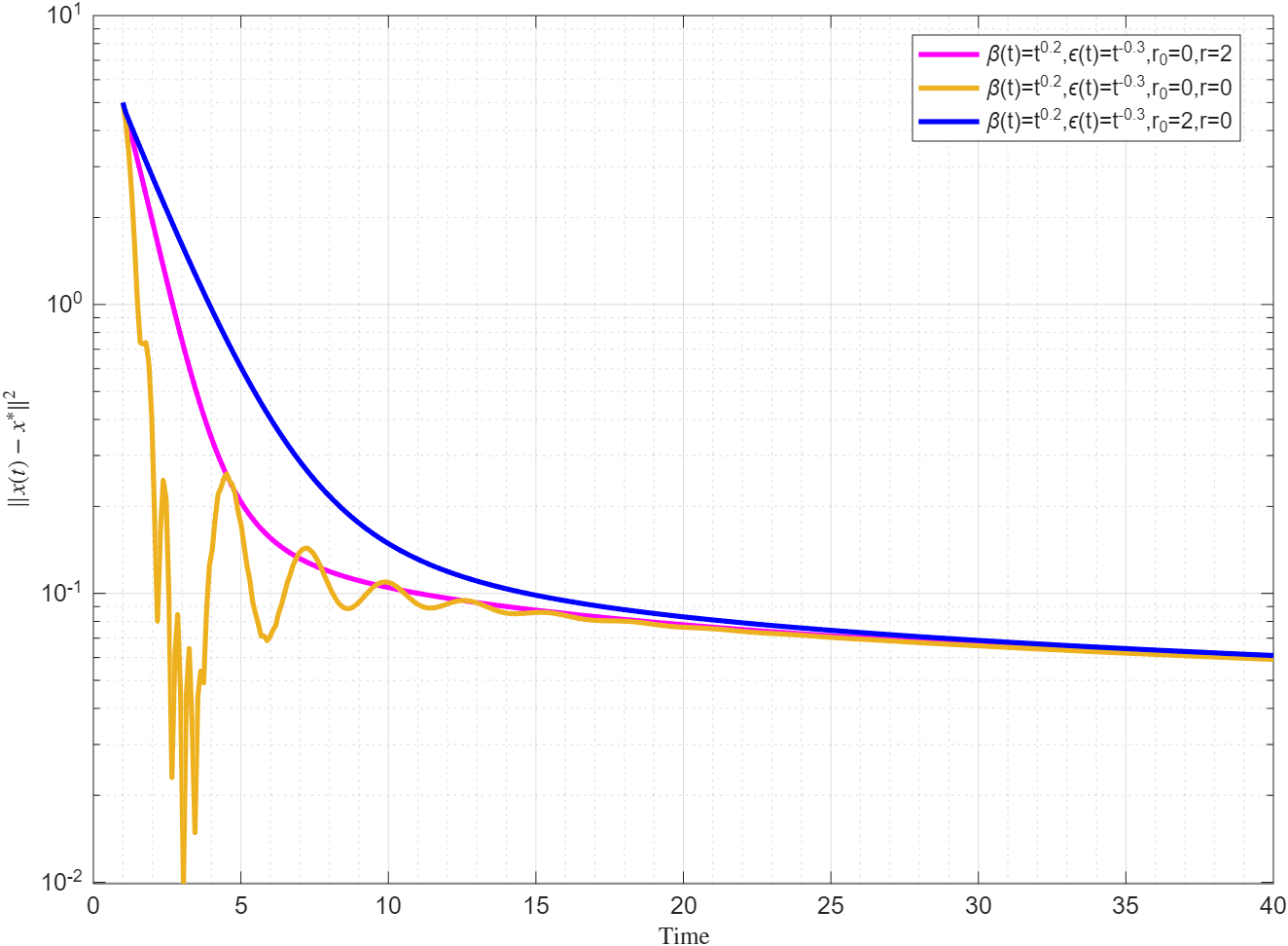}
    \caption{iterate error}
    \label{H:x}
\end{subfigure}
&
\begin{subfigure}{0.31\textwidth}
    \centering
    \includegraphics[width=\linewidth]{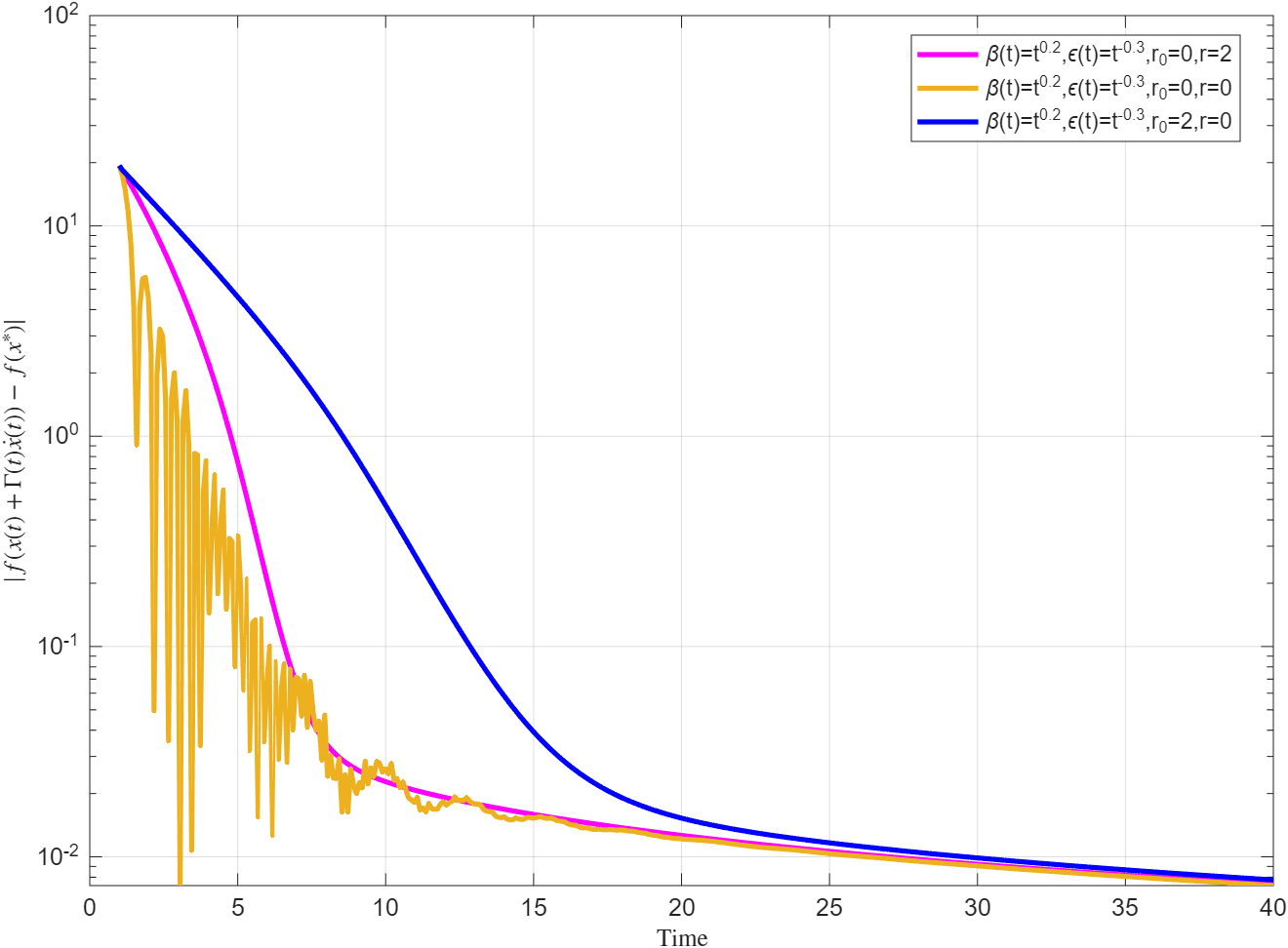}
    \caption{objective gap}
    \label{H:f}
\end{subfigure}
&
\begin{subfigure}{0.31\textwidth}
    \centering
    \includegraphics[width=\linewidth]{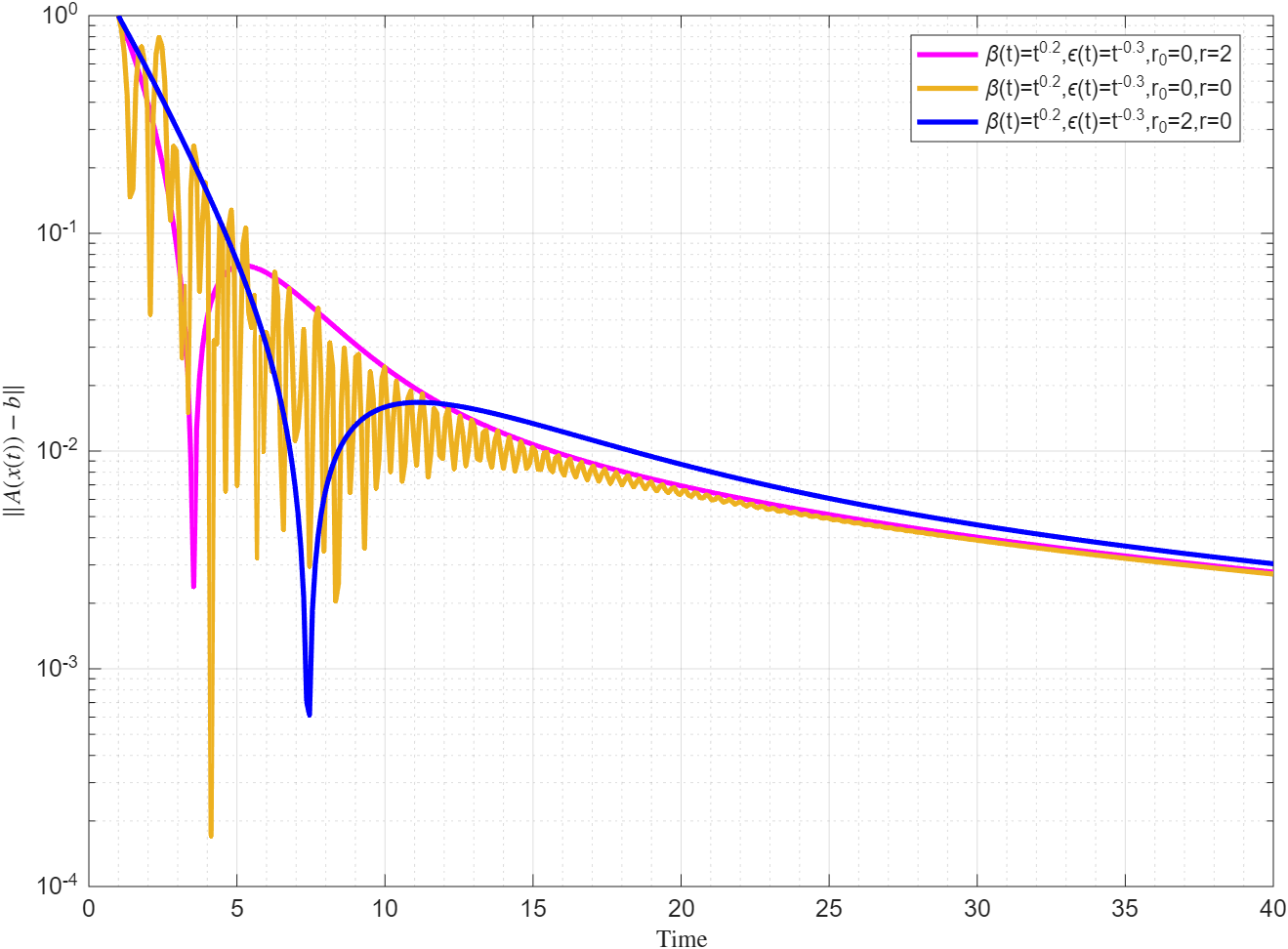}
    \caption{feasibility violation}
    \label{H:Ax-b}
\end{subfigure}
\end{tabular}

\caption{Comparison of the convergence behavior with and without implicit Hessian damping}
\label{fig:F4}
\end{figure}

As shown in Figure \ref{fig:F4}, when 
implicit Hessian damping term is removed, the feasibility violation, the objective residual, 
and the iterate error exhibit persistent oscillations throughout the evolution. In contrast, when the implicit Hessian damping mechanism is activated, the oscillations are 
substantially suppressed. These numerical results clearly demonstrate that the extrapolated gradient evaluation at 
$x(t)+\Gamma(t)\dot{x}(t)$ provides an effective implicit damping mechanism for reducing the 
oscillations.

\section{Conclusion}\label{conclusion}
In this paper, we proposed a mixed-order inertial primal-dual dynamical system with slowly vanishing damping, implicit Hessian damping, and the primal-dual Tikhonov regularization for solving linear equality-constrained convex optimization problems in Hilbert spaces. The proposed framework allows for general time-dependent time-scaling and Tikhonov regularization coefficients, rather than restricting them to prescribed power-law forms. Under suitable assumptions on the parameters and coefficient functions, a Lyapunov analysis yields convergence-rate estimates for the primal-dual gap and the objective residual at the extrapolated primal
point, as well as for the feasibility violation. We also establish the strong convergence of the generated primal-dual trajectories to the minimum-norm element of the saddle-point set. Finally, the numerical experiments illustrate the convergence behavior of the proposed dynamical system, demonstrate the oscillation-suppression effect of the implicit Hessian damping mechanism, and the minimum-norm selection effect induced by Tikhonov regularization.

\bibliographystyle{plain}

\bibliography{main}

\end{document}